\documentclass[a4paper]{article}
\usepackage[utf8]{inputenc}
\usepackage[T1]{fontenc}
\usepackage{lmodern}

\usepackage[a4paper]{geometry}
\usepackage{graphicx}
\usepackage[dvipsnames]{xcolor}
\usepackage{mathtools,amssymb}
\usepackage{todonotes}
\usepackage{hyperref}

\usepackage{algorithm}
\usepackage{algpseudocode}
\usepackage{colortbl}

\title{Multi-patient Inverse Estimation of Effective Membrane Diffusion Coefficients in Calcium–Citrate Hemodialysis}
\author{ \begin{minipage}[t]{0.45\textwidth}
\centering
\textbf{Geoffrey Lacour} \vspace*{0.2cm}\\
\small
Université Paris-Saclay\\
INRAE, MaIAGE\\
78350 Jouy-en-Josas, France\\
\texttt{geoffrey.lacour@inrae.fr}
\end{minipage}\and 
\begin{minipage}[t]{0.45\textwidth}
\centering
\textbf{Nicolae Cîndea} \vspace*{0.2cm}\\
\small
Université Clermont Auvergne,\\
CNRS, LMBP, F-63000 Clermont-Ferrand, France\\
\texttt{nicolae.cindea@uca.fr}
\end{minipage} \and 
\begin{minipage}[t]{0.45\textwidth}
\centering
\textbf{Julien Aniort} \vspace*{0.2cm} \\
\small
Département de Néphrologie, Dialyse et Transplantation,\\
CHU Clermont-Ferrand, Hôpital Gabriel Montpied,\\
58 Rue Montalembert, 63000, Clermont-Ferrand, France \vspace*{0.1cm}\\
\& \vspace*{0.1cm}\\
INRAE UMR 1019,\\
Unité de Nutrition Humaine, Université Clermont Auvergne,\\
Clermont-Ferrand, France\\
\texttt{janiort@chu-clermontferrand.fr}
\end{minipage}}

\begin{document}
\maketitle

\begin{abstract}
	We propose a multi-patient inverse modeling framework for identifying effective calcium and citrate diffusion coefficients in hollow-fiber hemodialysis devices. The approach relies on a coupled forward model combining axisymmetric fluid dynamics with multi-species convection-reaction-diffusion, together with a derivative-free optimization strategy to estimate membrane transport parameters from outlet concentration measurements. To account for inter-patient variability, physiological input parameters are first generated from clinical data and complemented by a patient-specific hydraulic calibration step, ensuring physical consistency across the synthetic cohort. The inverse problem is formulated as a global least-squares minimization aggregating residuals over multiple patients. Numerical experiments on synthetic data demonstrate multi-patient identifiability of the diffusion coefficients in the exact-data setting. Robustness with respect to measurement noise is subsequently assessed by perturbing observable outputs at various noise levels, and sensitivity analyses are performed to quantify the influence of membrane transport parameters on model predictions. The methodology is then applied to real clinical data obtained from an AK200 Gambro/Nikkiso DBB07 dialysis system. The results indicate that aggregating information from several patients substantially improves parameter identifiability and stability compared to single-patient inversions. Overall, this work provides a physically consistent and computationally tractable framework for multi-patient parameter estimation in dialysis models, and opens perspectives for large-scale personalization through physics-informed surrogate modeling.
\end{abstract}


\section{Introduction}
\label{sec:introduction}

Hemodialysis is a life-sustaining therapy for patients with end-stage renal disease (ESRD). Its primary objective is to restore physiological homeostasis by selectively removing or exchanging solutes and electrolytes between the blood and the dialysate across a semi-permeable membrane as described in J.Himmelfarb, T.Ikizler \& T.Alp~\cite{himmelfarb-ikizler-alp-2010}. Among these solutes, calcium plays a pivotal role in mineral metabolism and cardiovascular function. An inappropriate calcium balance during hemodialysis sessions may contribute to the development of chronic kidney disease–mineral and bone disorder (CKD-MBD), including osteomalacia, osteitis fibrosa, and adynamic bone disease. Moreover, dysregulated calcium homeostasis is strongly associated with cardiovascular complications such as valvular and vascular calcifications, as well as left ventricular hypertrophy (see M.Ketteler et al~\cite{ketteler-evenepoel-holden-isakova-jorgensen-2025}). 

In many clinical protocols, citrate is used as a regional anticoagulant owing to its capacity to chelate ionized calcium, thereby locally inhibiting activation of the coagulation cascade within the extracorporeal circuit. Accurate prediction of calcium transport and local concentration during intermittent hemodialysis therefore represents an important clinical and computational challenge. The coupling between citrate and calcium introduces additional complexity into transport dynamics, as calcium speciation depends on reversible binding reactions, convection, diffusion, and transmembrane exchanges. This coupling motivates detailed modeling of their interactions within the dialysis filter. Most commonly used models rely on one-dimensional approaches assuming equilibrium conditions. Although such models provide global estimates of calcium balance, they do not allow precise investigation of spatial variations in calcium concentration along the dialyzer fibers (see V.Maheshwari et al~\cite{maheshwari-cherif-fuertinger-schapachertilp-preciado-2017}). However, this spatial resolution is necessary to properly evaluate local anticoagulant efficacy, particularly since effective anticoagulation requires free calcium levels below 0.4 mmol/L within the filter. In previous work, we introduced a time- and space-dependent mathematical model of calcium transport in hollow-fiber dialyzers, based on coupled Navier–Stokes, Darcy, and convection–reaction–diffusion equations \cite{aniort-chupin-cindea-2018,aniort-richard-thouy-leguen-2024}. These models yield qualitatively realistic concentration profiles and are consistent with clinical observations. Nevertheless, their quantitative predictive capability strongly depends on membrane transport parameters, especially effective diffusion coefficients, which are difficult to measure directly and remain poorly characterized in practice.

In this work, we address the inverse problem of identifying effective diffusion coefficients in the dialysis membrane from sparse clinical measurements. More precisely, given inlet concentrations and flow rates, together with outlet blood concentrations for several chemical species, we aim at estimating membrane diffusion parameters so that numerical simulations best reproduce observed patient data. This constitutes a nonlinear, low-dimensional inverse problem driven by an expensive forward model, involving coupled fluid dynamics and convection-reaction-diffusion equations. Additional difficulties arise from limited observability, measurement noise, and interpatient variability. Beyond model calibration, identifying effective membrane transport coefficients is a key step toward patient-specific dialysis strategies. Reliable parameter estimates would enable predictive simulations under varying dialysate compositions, opening the possibility of optimizing treatment protocols to achieve targeted calcium balance while minimizing adverse effects. Such personalized approaches are increasingly sought in clinical practice, but require quantitative models grounded in measurable patient data.

The present study was conducted within the framework of the MARC clinical trial (ClinicalTrials.gov identifier: NCT04530175), carried out at CHU Gabriel Montpied (Clermont-Ferrand, France). The optimization of the mathematical model, was prespecified as a secondary objective of the study protocol. The analysis is based on data from twenty-two patients undergoing intermittent hemodialysis. Building upon the existing forward model, we formulate a multi-patient inverse problem and adopt a derivative-free optimization strategy to estimate two key diffusion coefficients associated with calcium and citrate transport across the membrane. To assess the reliability of the proposed approach, we first validate the method on synthetic data, allowing us to investigate identifiability in the noise-free case and robustness with respect to measurement perturbations. We then apply the methodology to clinical data obtained on AK200 Gambro and Nikkiso DBB07 dialysis machines using hollow-fiber dialyzer cartridges equipped with synthetic polysulfone membranes (Leoceed-21 and/or APS series, distributed by Asahi Kasei Medical Co., Ltd.), and discuss the resulting parameter estimates. 

The remainder of the paper is organized as follows. Section~\ref{sec:math-model} introduces the forward model, combining fluid dynamics with multi-species transport and chemical reactions in an axisymmetric hollow-fiber geometry, and details the numerical resolution of the stationary problem. Section~\ref{sec:inverse-pipeline} presents the formulation of the inverse problem together with its numerical implementation. In the latter, we first consider in Section~\ref{sec:single-patient} a single-patient inverse
problem as a baseline validation of the modeling and optimization pipeline. Then, Sections~\ref{subsec:synthetic-identifiability}--\ref{subsec:synthetic-sensitivity} address the multi-patient framework using synthetic data, focusing successively on exact-data identifiability, robustness with respect to measurement noise, and sensitivity to the biomembrane diffusion coefficients. Finally, Section~\ref{sec:real-data-inverse} presents an application to clinical data obtained dialysis devices mentioned above. In Section~\ref{sec:conclusions}, we end the presentation by summarizing the main findings and discuss perspectives, including the integration of physics-informed neural networks to reduce computational cost and improve scalability with respect to patient cohorts.

Although the present contribution builds upon the framework introduced in the works of the present authors \& L.Chupin~\cite{aniort-chupin-cindea-2018} as well as the work of the present authors with F.Richard, F.Thouy, L.Le Guen, C.Philipponnet, C.Garrouste, A.-E.Heng, C.Dupuis, M.Adda, D.Julie, E.Lebredonchel, L.Chupin, D.Bouvier \& B.Souweine~\cite{aniort-richard-thouy-leguen-2024}, it should also be viewed within a broader and highly active research landscape. Over the past decades, extensive efforts have been devoted to the modeling and the study of transport phenomena in membrane-based systems, to the mathematical analysis of related mathematical models, and to the identification of effective transport parameters from experimental or clinical data. Significant advances have been achieved in the areas of membrane characterization, inverse modeling, and parameter estimation in complex multiphysics settings, as illustrated for instance in the works of N.Cancilla, L.Gurreri, M.Ciofalo, A.Cipollina, A.Tamburini \& M.Giorgio~\cite{cancilla-gurreri-ciofalo-cipollina-tamburini-2025}, of N.Cancilla, L.Gurreri, G.Marotta, M.Ciofalo, A.Cipollina, A.Tamburini \& G.Micale~\cite{cancilla-gurreri-marotta-ciofalo-cipollina-2022}, of Y.-L.Chang \& C.-J.Lee~\cite{chang-lee-1988}, of D.Donato, A.Boschetti-de-Fierro, C.Zweigart, M.Kolb, S.Eloot, M.Storr, B.Krause, K.Leypoldt \& P.Segers~\cite{donato-boschetti-zweigart-kolb-eloot-2017}, of D.Donato, M.Storr \& B.Krause~\cite{donato-storr-krause-2020}, of S.Eloot, J.Vierendeels \& P.Verdonck~\cite{eloot-vierendeels-verdonck-2006}, of S.Eloot, D.De Wachter, I.Van Tricht \& P.Verdonck~\cite{eloot-dewatcher-vantricht-verdonck-2002}, of A.Fasano \& A.Farina~\cite{fasano-farina-2012}, of M.Islam \& J.Szpunar~\cite{islam-szpunar-2013}, of M.Jaffrin, B.Gupta \& J.Malbrancq~\cite{jaffin-gupta-malbrancq-1981}, of C.Legallais, G.Catapano, B.von Harten \& U.Baurmeister~\cite{legallais-catapano-vonharten-baurmeister-2000}, of T.Miyasaka \& K.Sakai~\cite{miyasaka-sakai-2023}, of C.Ronco \& W.Clark~\cite{ronco-clark-2018}, of D.Schneditz, D.Platzer \& J.Daugirdas~\cite{schneditz-platzer-daugirdas-2009}, of A.Werynski \& J.Waniewski~\cite{werynski-waniewski-1995}, of T.Yaqoob, M.Ahsan, A.Hussain \& I.Ahmad~\cite{yaqoob-ahsan-hussain-ahmad-2020}, of J.Zhang, X.Chen, J.Ding, K.Fraser, M.Ertan Taskin, B.Griffith \& Z.Wu~\cite{zhang-chen-ding-faser-ertantaskin-2013} and references therein. We may also refer to the book of R.Fournier~\cite{fournier-2011} as well as the work of J.Himmelfarb, T.Ikizler \ \& T.Alp~\cite{himmelfarb-ikizler-alp-2010}. The present work contributes to this ongoing line of research by combining mechanistic partial differential equation modeling, multi-patient inverse identification, and numerical sensitivity analysis in the specific context of hollow-fiber hemodialysis devices.

\section{Forward model: fluid dynamics and multi-species transport in a hollow-fiber dialyzer}\label{sec:math-model}

This section presents the forward model that underlies the inverse identification procedure developed later. Our approach follows the modeling framework introduced in the work of the present authors \& their collaborators~\cite{aniort-chupin-cindea-2018, aniort-richard-thouy-leguen-2024}. We recall the main assumptions, describe the governing equations in an axisymmetric setting, and explain the boundary conditions and numerical strategy used to compute the outlet concentrations. For the sake of clarity and reproducibility, we emphasize the physical meaning of each model component, while postponing the full non-dimensionalization details and parameter tables to \cite{aniort-chupin-cindea-2018}.

\begin{figure}[H]
  \centering
  \includegraphics[width=0.5\textwidth]{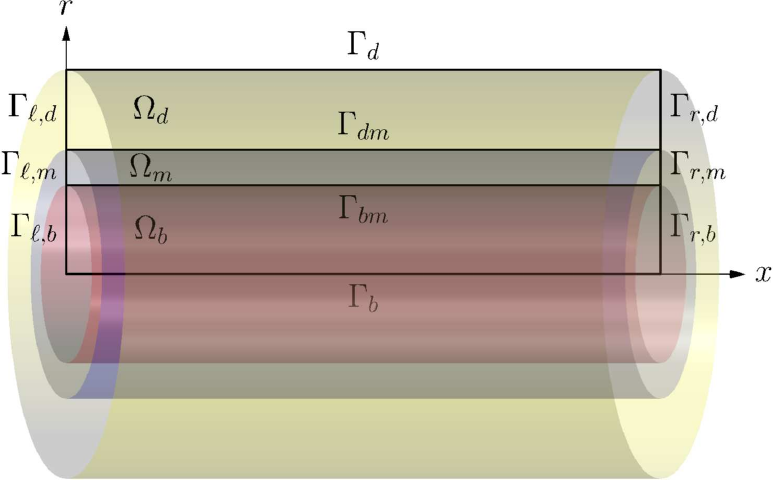}
  \caption{Axisymmetric representation of a representative hollow fiber: blood region $\Omega_b$, porous membrane $\Omega_m$, and dialysate region $\Omega_d$.}
  \label{fig:omega}
\end{figure}

A dialyzer module contains a large number of parallel hollow fibers. A classical modeling approach consists in studying a representative fiber and assuming periodicity of the surrounding structure. Each fiber may be divided into three regions: the blood channel (inner lumen), a porous membrane and the dialysate channel (outer region). We may consider as a suitable approximation that the fiber geometry is cylindrical so that flows are approximately radially symmetric, letting us work in cylindrical coordinates and reduce the three-dimensional setting to a two-dimensional axisymmetric domain in variables $(x,r)$. Having this in mind, we denote by $L$ the fiber length and by $R$ its outer radius. Let $R_1$ be the blood channel radius and $R_2$ the outer membrane radius.The computational domain is the rectangle
\begin{equation*}
	\Omega = \{(x,r)\in\mathbb{R}^2\, |\, \ 0<x<L,\ 0<r<R\},
\end{equation*}
and we can decompose it into
\begin{equation*}
\Omega_b=\{0<r<R_1\},\qquad \Omega_m=\{R_1<r<R_2\},\qquad \Omega_d=\{R_2<r<R\},
\end{equation*}
so interfaces are given by $\Gamma_{bm}=\{r=R_1\}$ and $\Gamma_{dm}=\{r=R_2\}$. Also, the boundary is decomposed into inlet/outlet parts for each subdomain and two horizontal boundaries corresponding to the symmetry axis $\{r=0\}$ and the outer confinement boundary $\{r=R\}$ (see Figure~\ref{fig:omega}).

Using an axisymmetric formulation in cylindrical coordinates properly accounts for the fiber geometry. In particular, conservation laws naturally lead to operators of the form $\partial_x +\frac{1}{r}\partial_r(r\,\cdot)$, while cross-sectional quantities involve the radial weight $r$. This is essential to compute physically meaningful outlet concentrations and to enforce regularity and symmetry conditions at the axis $r=0$.

\subsection{Flow Model and Velocity Field Computation}
\label{subsec:flow}

Solute transport is driven by the velocity field in the blood and dialysate compartments. Because hollow fibers are well approximated by cylindrical geometries, we formulate the flow model in cylindrical coordinates $(x,r,\theta)$, where $x$ denotes the axial direction, $r$ the radial variable, and $\theta$ the azimuth. We denote by $\mathbf v=(v_x,v_r,v_\theta)$ the fluid velocity and by $p$ the pressure. As a reference continuum model, we consider the incompressible Navier-Stokes equations for a Newtonian fluid, which in cylindrical coordinates are given by
\begin{equation*}
	\begin{cases}
		\begin{split}
			\rho\bigg(\partial_t v_x + v_x\partial_x v_x + v_r\partial_r v_x + \frac{v_\theta}{r}\partial_\theta v_x\bigg) =&\, -\partial_x p + \mu\bigg(\partial_x^2 v_x+\frac{1}{r}\partial_r(r\partial_r v_x)+\frac{1}{r^2}\partial_\theta^2 v_x\bigg),\vspace*{0.2cm}\\
			\rho\bigg(\partial_t v_r + v_x\partial_x v_r + v_r\partial_r v_r +\frac{v_\theta}{r}\partial_\theta v_r - \frac{v_\theta^2}{r}\bigg) = &\, \mu\bigg(\partial_x^2 v_r+\frac{1}{r}\partial_r(r\partial_r v_r) + \frac{1}{r^2}\partial_\theta^2 v_r-\frac{v_r}{r^2}-\frac{2}{r^2}\partial_\theta v_\theta\bigg)\vspace*{0.1cm}\\
			& \, -\partial_r p, \vspace*{0.2cm} \\
			\rho\bigg(\partial_t v_\theta + v_x\partial_x v_\theta + v_r\partial_r v_\theta + \frac{v_\theta}{r}\partial_\theta v_\theta + \frac{v_r v_\theta}{r}\bigg) =&\,\mu\bigg(\partial_x^2 v_\theta+ \frac{1}{r}\partial_r(r\partial_r v_\theta)+\frac{1}{r^2}\partial_\theta^2 v_\theta \frac{v_\theta}{r^2}+\frac{2}{r^2}\partial_\theta v_r\bigg)\vspace*{0.1cm}\\
			&\,  -\frac{1}{r}\partial_\theta p, \vspace*{0.2cm}\\
			\partial_x v_x+\frac{1}{r}\partial_r(r v_r)+\frac{1}{r}\partial_\theta v_\theta=&\,0,
		\end{split}
	\end{cases}
\end{equation*}
where $\rho$ is the (constant) density and $\mu$ the dynamic viscosity.
We refer to standard textbooks in fluid mechanics for the derivation of these expressions, see e.g. the book of P. Davidson~\cite{davidson-2024}. In hollow-fiber dialyzers, both the geometry and operating conditions are well approximated by an axisymmetric setting, and no significant azimuthal motion is observed. Mathematically, this is represented through $\partial_\theta(\cdot)=0$ and $v_\theta\equiv0$. Under such assumptions, the velocity reduces to $\mathbf v=(v_x(x,r),v_r(x,r))$ and the incompressibility condition, which is the ultimate equation in the above Navier-Stokes system, reduces to
\begin{equation*}
	\partial_x v_x+\frac1r\partial_r(r v_r)=0.
\end{equation*}
The momentum equations simplify accordingly, yielding the classical axisymmetric no-swirl Navier-Stokes system
\begin{equation*}
	\begin{cases}
		\begin{split}
			\rho\bigg(\partial_t v_x + v_x\partial_x v_x + v_r\partial_r v_x\bigg) =&\, -\partial_x p + \mu\bigg(\partial_x^2 v_x+\frac{1}{r}\partial_r(r\partial_r v_x)\bigg),\vspace*{0.2cm}\\
			\rho\bigg(\partial_t v_r + v_x\partial_x v_r + v_r\partial_r v_r\bigg) = &\, -\partial_r p + \mu\bigg(\partial_x^2 v_r + \frac{1}{r}\partial_r(r\partial_r v_r) -\frac{v_r}{r^2}\bigg)\vspace*{0.2cm}\\
			\partial_x v_x+\frac{1}{r}\partial_r(r v_r)=&\,0.
		\end{split}
	\end{cases}
\end{equation*}
At the scale of a single fiber, characteristic radii are of order $R \sim 10^{-4}\,$m, while module-level blood and dialysate flow rates are distributed over approximately $N \sim 10^4$ parallel fibers. Using typical clinical values reported in the work of the present authors \& L.Chupin \cite{aniort-chupin-cindea-2018}, this leads to moderate Reynolds numbers at the fiber level, consistent with a laminar flow regime. In the present work, we therefore neglect inertial effects and solve a stationary Stokes-type creeping flow problem in both the blood and dialysate compartments. This may be surprising at fist sight, since blood may exhibit non-Newtonian behavior in general. However, in narrow conduits the apparent viscosity decreases due to the F{\aa}hr{\ae}us-Lindqvist effect, and the flow can therefore be accurately approximated by a Newtonian model. Moreover, numerical comparisons reported in \cite{aniort-chupin-cindea-2018} indicate that predicted solute concentrations are only weakly sensitive to the specific rheological law in the regimes considered. This motivates the use of a Newtonian Stokes model for both blood and dialysate in the inverse identification pipeline.

Thus, the above assumptions imply that both the time derivatives and the convective-type terms can be neglected, since we consider an axisymmetric and laminar regime without swirl, so that the governing equations reduce to
\begin{equation}\label{eq:fluid-flow}\tag{F$_c$}
	\begin{cases}
		\mu\bigg(\partial_x^2 v_x+\frac{1}{r}\partial_r(r\partial_r v_x)\bigg) -\partial_x p  = 0, \vspace*{0.1cm}\\
		\mu\bigg(\partial_x^2 v_r+\frac{1}{r}\partial_r(r\partial_r v_r)-\frac{v_r}{r^2}\bigg) -\partial_r p = 0, \vspace*{0.1cm}\\
		\partial_x v_x+\frac{1}{r}\partial_r(r v_r)=0,
	\end{cases}
\end{equation}
posed in the blood and dialysate domains with accordiing parameters. The hollow-fiber membrane is modeled as a porous medium. In this region, inertial effects remain negligible and the flow is assumed to obey Darcy law. Denoting by $p_m$ the pressure in the membrane and by $\mathbf u_m=(u_{m,x},u_{m,r})$ the fluid velocity accross the membrane, the governing equations can be written in axisymmetric cylindrical coordinates
\begin{equation}\label{eq:membrane-darcy}\tag{F$_m$}
	\begin{cases}
		u_{m,x}=-\frac{K}{\mu}\partial_x p_m, \vspace*{0.1cm}\\
		u_{m,r}=-\frac{K}{\mu}\partial_r p_m, \vspace*{0.1cm} \\
		\partial_x u_{m,x}+\frac{1}{r}\partial_r(r u_{m,r})=0,
	\end{cases}
\end{equation}
so that the membrane pressure does satisfy the elliptic equation
\begin{equation}\label{eq:pressure-membrane}\tag{P$_m$}
	\partial_x^2 p_m + \frac{1}{r}\partial_r(r\partial_r p_m)=0 \quad \text{in}\; \Omega_m.
\end{equation}
Physically, Darcy law expresses a linear relation between pressure gradients and filtration velocity in the porous membrane, with the permeability $K$ encoding the microstructural properties of the fiber wall. This model is appropriate for slow flows in highly resistive porous media and is usual in dialyzer-scale modeling. The Stokes systems in $\Omega_b$ and $\Omega_d$ given by \eqref{eq:fluid-flow} are coupled to the Darcy flow in $\Omega_m$
given in \eqref{eq:membrane-darcy} and \eqref{eq:pressure-membrane} through interface conditions enforcing continuity of the normal flux and appropriate stress or pressure transmission. In practice, fiber-level inlet and outlet pressures are not directly measured, whereas blood and dialysate flow rates are available at the machine level. Following \cite{aniort-chupin-cindea-2018}, we determine consistent pressure boundary values by solving a calibration problem so that the simulated fluxes match the measured ones. The resulting velocity field is then kept fixed and used as input for the solute transport model described in the next subsection. We refer to \cite{aniort-chupin-cindea-2018} for further implementation details.

\subsection{Multi-species Transport and Chemical Reactions }
\label{subsec:transport}

We consider the transport and interaction of five chemical species in blood and dialysate that play a central role in calcium balance under citrate anticoagulation, as described in \cite{aniort-chupin-cindea-2018}: free calcium ($c_1$), free albumin binding sites ($c_2$), calcium-albumin complexes ($c_3$), citrate ($c_4$), and calcium-citrate complexes ($c_5$). These species are known to govern ionized calcium availability and buffering mechanisms during dialysis sessions under citrate anticoagulation, see e.g. the textbook of J.Hall \& M.Hall~\cite[Unit V]{hall-hall-2020} for a clinical overview.

We denote by $\mathbf c=(c_1,c_2,c_3,c_4,c_5)^T$ the vector of concentrations.
Then for each species, conservation of mass over a control volume yields a balance equation of the form
\begin{equation*}
	\partial_t c_i + \mathrm{div}\big(\mathbf J_i\big) = F_i(\mathbf c), \quad (i=1,\ldots,5)
\end{equation*}
where $\mathbf J_i$ is the total flux of species $i$ and $F_i(\mathbf c)$ represents the net production rate due to local chemical reactions. Namely, it accounts for reversible binding reactions between calcium, albumin and citrate, modeled through mass-action kinetics. Such reaction-diffusion systems naturally lead to coupled convection-reaction-diffusion equations and are classical in the modeling of reactive solute transport in biological fluids (see, e.g. G.Truskey, F.Yuan \& D.Katz~\cite{truskey-yuan-katz-2008} for a general modeling perspective. Referring to \cite{aniort-chupin-cindea-2018} the nonlinear source term which account for interaction between chemical species can is written
\begin{equation}\label{eq:interacting-term}
  F(\boldsymbol{c}) = \begin{pmatrix}F_1(\boldsymbol{c})\\ F_2(\boldsymbol{c})\\ F_3(\boldsymbol{c})\\ F_4(\boldsymbol{c}) \\ F_5(\boldsymbol{c})\end{pmatrix}= \frac{1}{\mathcal{F}d}
  \begin{pmatrix}
    c_3 + \delta_2 c_5 - \delta_1 c_1 c_2 - \delta_3 c_1 c_4 \\
    c_3 - \delta_1 c_1 c_2 \\
    - c_3 + \delta_1 c_1 c_2 \\
    \delta_2 c_5 - \delta_3 c_1 c_4 \\
    -\delta_2 c_5 + \delta_3 c_1 c_4
  \end{pmatrix},
\end{equation}
where $\delta_1$, $\delta_2$, $\delta_3$ and $\mathcal{F}d$ are some constants. More precisely, the reaction network corresponds to reversible bindings
\begin{equation*}
	\text{Ca}+\text{Alb} \rightleftharpoons \text{Ca-Alb} \quad \text{and} \quad \text{Ca}+\text{Cit} \rightleftharpoons \text{Ca-Cit},
\end{equation*}
modeled through first-order kinetics with constants that can be associated here (up to a rescaling method) to $\delta_1, \ldots, \delta_3$. Let us focus on the flux term $\mathbf J_i$. We decompose it as the sum of an advective flux and a diffusive flux:
\begin{equation*}
	\mathbf J_i = S_i\,c_i\,\mathbf U - D_i \nabla c_i.
\end{equation*}
Here $\mathbf U = (U_x,U_r)$ is the velocity field solution to \eqref{eq:fluid-flow}. The coefficient $S_i$ is a sieving factor (which is dimensionless), used to represent hindered transport across the membrane: $S_i=1$ corresponds to unhindered solute transport, whereas $S_i=0$ suppresses convection in the membrane for large molecules such as albumin. The diffusion coefficient $D_i$ may take different values across $\Omega_b$, $\Omega_m$, and $\Omega_d$. In cylindrical coordinates for an axisymmetric setting, such a consideration of $\mathbf J_i$ yields the convection-reaction-diffusion system
\begin{equation*}
	\partial_t c_i + S_i\,(U_x\partial_x c_i + U_r\partial_r c_i) - \partial_x(D_i\partial_x c_i) -\frac{1}{r}\partial_r\!\left(r D_i \partial_rc_i\right) = F_i(\mathbf c), \qquad (i=1,\ldots,5)
\end{equation*}
posed on the whole domain $\Omega$. Although $D_i$ and $S_i$ may be discontinuous across subdomains, all quantities are defined on $\Omega$ and the variational formulation remains meaningful at the global level, as described in  \cite{aniort-chupin-cindea-2018}. 

To benefit from a suitable description of the physical underlying phenomenon, we do need to impose boundary conditions that reflect known inlet concentrations and the fact that there is no loss of chemical species in the whole device. To this end, we assume that blood concentrations at the inlet are measured and dialysate concentrations at the inlet are known from the dialysate prescription. This yields Dirichlet boundary conditions given by
\begin{equation*}
	c_i = \overline{c_{i,\text{input}}} \ \text{ on } \Gamma_{\ell,b} \cup \Gamma_{r,d} \qquad (i=1,\ldots,5).
\end{equation*}
Such conditions represent the measured (or controlled) chemical state into the model. On the remaining boundaries, concentrations are not prescribed. A common modeling choice is to assume that there is no loss through the exterior boundary of the device, which yields homogeneous Neumann conditions
\begin{equation*}
\partial_{\boldsymbol n} c_i = 0 \ \text{ on } \Gamma_{b} \cup \Gamma_{r,b} \cup \Gamma_{\ell, d} \quad (i = 1,4,5) \quad \text{and} \quad \partial_{\boldsymbol n} c_j = 0 \ \text{ on } \Gamma_{b} \cup \Gamma_{r,b} \quad (j = 2,3).
\end{equation*}
Together, these conditions mix prescribed values (Dirichlet) where the chemical composition is known, and null flux assumptions (Neumann) where it is not. This mixed structure is essential for identifiability since the model is driven by inputs, while outputs remain free and become meaningful observables to compare with clinical measurements. However, it may yield several mathematical issues both from a computational and a theoretical point of views.

Finally, let us remember that we consider ourselves to be in an equilibrium state, which means that we can consider the exchange of chemical species in the dialyzer intrinsically, and thus view the system as being in equilibrium (without time derivative). Namely, this means that the regimes considered, the transient system relaxes rapidly toward a steady state along the fiber. Gathering the above observations, we may describe the  considered system to describe chemical interactions into a hollow-fiber of the dialyzer:
\begin{itemize}
\item For $i \in \{1, 4, 5\}$ which we recall respectively corresponds to Calcium, Citrate and Calcium-Citrate, the non-transient system is given as
\begin{equation}\label{eq:dar145}
	\begin{cases}
    		U_x \partial_x c_i + U_r \partial_r c_i  - \frac{1}{\mathcal{P}e} \frac{1}{r} \partial_r
   		\left(
    rD_i \partial_r c_i
    \right)
      - \frac{\varepsilon^2}{\mathcal{P}e} \partial_x \left(D_i \partial_x c_i \right) = F_i(\boldsymbol{c})& \text{ in } \Omega \\
    		c_i = \overline{c_{i,\ \text{input}}}& \text{ on }  \Gamma_{\ell,b} \cup \Gamma_{r, d}\\
    		\partial_{\boldsymbol{n}}c_i = 0& \text{ on } \Gamma_{b} \cup \Gamma_{r,b} \cup \Gamma_{\ell, d};

  \end{cases}
\end{equation}
  
\item Also since Albumin and Calcium-Albumin do not cross the membrane, the equations representing species $i \in \{2,\ 3\}$ are fully described in $\Omega_b$ and are given by
\begin{equation}\label{eq:dar23}
  \begin{cases}
    U_x \partial_x c_i + U_r \partial_r c_i  -
    \frac{1}{\mathcal{P}e}
    \frac{1}{r} \partial_r
    \left(
    rD_i \partial_r c_i
    \right)
    -
    \frac{\varepsilon^2}{\mathcal{P}e}
    \partial_x \left(D_i \partial_x c_i \right) = F_i(\boldsymbol{c})& \text{ in } \Omega_b \\
    c_i = \overline{c_{i,\ \text{input}}}& \text{ on }  \Gamma_{\ell,b} \\
    \partial_{\boldsymbol{n}}c_i = 0& \text{ on } \Gamma_{b} \cup \Gamma_{r,b}\\
    D_i \partial_{r}c_i = c_i U_r& \text{ on } \Gamma_{bm}.
  \end{cases}
\end{equation}
\end{itemize}
In systems~\eqref{eq:dar145} and \eqref{eq:dar23}, the interaction terms are given by \eqref{eq:interacting-term}, and all constants correspond to physical parameters (e.g. a Péclet number and a $\varepsilon^2$-weighted longitudinal diffusion term). Since their explicit values are not needed for the present study, we refer to \cite{aniort-chupin-cindea-2018} for a detailed description.

\subsection{Numerical Resolution of the Forward Problem}
\label{subsec:forward-stationary}

For each species indexed by $i$, the predicted outlet concentration is computed f or each species using the formula  
\begin{equation}\label{eq:output-formmula}\tag{OP}
	c_{i,\ \text{output}} = \frac{2R^2}{R_1^2} \int_{\Gamma_{r, b}} r c_i dr,
\end{equation}
namely the model output is defined as the cross-sectional average concentration in blood at $x=L$ for a solution $\boldsymbol{c}$ of \eqref{eq:dar145}--\eqref{eq:dar23}. We discretize $\Omega$ with a structured triangulation and approximate the concentrations using $P_1$ finite elements. Note that the stationary coupled system \eqref{eq:dar145}--\eqref{eq:dar23} is nonlinear due to the reaction term $F(\mathbf c)$. We approach the solution of the stationnary convection-reaction-diffusion system using a Newton method in a variational form. More exactly we first write the variational formulation associated to~\eqref{eq:dar145}, so that for $i \in \{1,\ 4, \ 5\}$, we aim to find $c_i \in  \mathcal{C}_i$ solution of
\begin{align}
\iint_\Omega \left( U_x \partial_x c_i + U_r \partial_r c_i \right) r \varphi_i \, \mathrm dx \, \mathrm dr + \frac{1}{\mathcal{P}e} \iint_\Omega r D_i \partial_r c_i \partial_r \varphi_i \, \mathrm dx \, \mathrm dr - \frac{\varepsilon^2}{\mathcal{P}e} \int_{\Gamma_{r,d}} r D_i \partial_x c_i \varphi_i \, \mathrm dr \nonumber \\
 + \frac{\varepsilon^2}{\mathcal{P}e} \int_{\Gamma_{\ell, b}} r D_i \partial_x c_i \varphi_i \, \mathrm dr + \frac{\varepsilon^2}{\mathcal{P}e} \iint_\Omega r D_i \partial_x c_i \partial_x \varphi_i \, \mathrm dx \, \mathrm dr = \iint_\Omega  r F_i(\boldsymbol{c}) \varphi_i \, \mathrm dx \, \mathrm dr, \label{eq:vf145}
\end{align}
for every $\varphi_i \in \mathcal{C}_i$, where $\mathcal{C}_i$ are suitable Sobolev spaces, that in particular take into account the boundary conditions described in~\eqref{eq:dar145}. Similarly, we can write a variational formulation associated to~\eqref{eq:dar23}. The main difference with respect to the variational formulation~\eqref{eq:vf145} is that the solutions to~\eqref{eq:dar23} are supported in $\Omega_b$ and not in $\Omega$ and are null in $\Omega_m$ and $\Omega_d$. Another difference is the presence of a boundary term on $\Gamma_{bm}$ due to the boundary condition. More precisely, for $i \in \{2, 3\}$ are looking for $c_i \in \mathcal{C}_i$ which solve
\begin{align}
  \iint_{\Omega_b} \left( U_x \partial_x c_i + U_r \partial_r c_i \right) r \varphi_i \, \mathrm dx \, \mathrm dr + \frac{1}{\mathcal{P}e} \iint_{\Omega_b} r D_i \partial_r c_i \partial_r \varphi_i \, \mathrm dx \, \mathrm dr
- \int_{\Gamma_{bm}} \frac{U_r}{\mathcal{P}e}c_i \varphi_i \mathrm dr
  \nonumber \\
 + \frac{\varepsilon^2}{\mathcal{P}e} \int_{\Gamma_{\ell, b}} r D_i \partial_x c_i \varphi_i \, \mathrm dr + \frac{\varepsilon^2}{\mathcal{P}e} \iint_{\Omega_b} r D_i \partial_x c_i \partial_x \varphi_i \, \mathrm dx \, \mathrm dr = \iint_{\Omega_b}  r F_i(\boldsymbol{c}) \varphi_i \, \mathrm dx \, \mathrm dr, \label{eq:vf23}
\end{align}
for every $\varphi_i \in \mathcal{C}_i$, where once again $\mathcal{C}_i$ are Sobolev spaces taking into account the boundary conditions present in~\eqref{eq:dar23}. We highlight that we can gather \eqref{eq:vf145}-\eqref{eq:vf23}, and rewrite the variational problem so that we are seeking for $\boldsymbol{c} = (c_i)_{1 \le i \le 5} \in \boldsymbol{\mathcal{C}}$ solving
\begin{equation}
  \label{eq:vf}
  \mathcal{F}(\boldsymbol{c}, \boldsymbol{\varphi}) = \boldsymbol{0},
\end{equation}
for every $\boldsymbol{\varphi} = (\varphi_i)_{1 \le i \le 5} \in \boldsymbol{\mathcal{C}}=\mathop{\times}\limits_{i=1}^5 \mathcal{C}_i$ with $\mathcal{F} : \boldsymbol{\mathcal{C}} \times \boldsymbol{\mathcal{C}} \to \mathbb{R}^5$ regrouping the integral terms in the five equalities~\eqref{eq:vf145}-\eqref{eq:vf23}.

To implement a Newton method, we have to compute the gradient of $\mathcal{F}$, denoted as usual $\nabla \mathcal{F}$. Moreover, we can write it as follows
\begin{equation}\label{eq:gradF}
\nabla \mathcal{F}(\boldsymbol{c}, \boldsymbol{\varphi}) \boldsymbol{\overline{c}} = \lim_{h \to 0} \frac{\mathcal{F}(\boldsymbol{c} + h\boldsymbol{\overline{c}}, \boldsymbol{\varphi}) - \mathcal{F}(\boldsymbol{c}, \boldsymbol{\varphi})}{h} = \mathcal{L}(\boldsymbol{\overline c}, \boldsymbol{\varphi}) + \mathcal{N}(\boldsymbol{c}, \boldsymbol{\overline c}, \boldsymbol{\varphi}).
\end{equation}
In the above, the term $\mathcal{L}$ brings together the linear terms in the right-hand sides of~\eqref{eq:dar145}--\eqref{eq:dar23}, while $\mathcal{N}$ gathers the remaining nonlinear terms, \emph{i.e.},
\[
\mathcal{N}(\boldsymbol{c}, \boldsymbol{\overline c}, \boldsymbol{\varphi}) = \left( \iint_\Omega r (\nabla F(\boldsymbol{c})\boldsymbol{\overline c})_i \varphi_i \, \mathrm dx \, \mathrm dr \right)_{1 \le i \le 5}.
\]

We give the detail of the variational Newton method employed to approach the solution of~\eqref{eq:vf} in Algorithm~\ref{alg:newton}. 
\begin{algorithm}[H]
\caption{Newton method in variational form for equation~\eqref{eq:vf}.}\label{alg:newton}
\begin{algorithmic}
\Require $\mathcal{F}$, $\boldsymbol{c_0}$, $tol$, $n_\text{max}$ \Comment{$n_\text{max}$ is the maximal number of iterations}
\State $n \gets 0$.
\State Compute $\boldsymbol{c}_{n+1} \in \boldsymbol{\mathcal{C}}$ solution to
\[
  \nabla \mathcal{F}(\boldsymbol{c_n}, \boldsymbol{\varphi}) \boldsymbol{c_{n+1}} = \nabla \mathcal{F}(\boldsymbol{c_n}, \boldsymbol{\varphi}) \boldsymbol{c_n} - \mathcal{F}(\boldsymbol{c_n}, \boldsymbol{\varphi}) \qquad \text{ for every } \boldsymbol{\varphi} \in \boldsymbol{\mathcal{C}}.
\]
\While{$n \le n_\text{max}$  and $\|\boldsymbol{c_{n+1}} - \boldsymbol{c_{n}}\| > \text{tol}$} \Comment{The tolerance tol is taken equal to $10^{-4}$} 
\State $n \gets n + 1$
\State Compute $\boldsymbol{c}_{n+1} \in \boldsymbol{\mathcal{C}}$ solution to
\[
  \nabla \mathcal{F}(\boldsymbol{c_n}, \boldsymbol{\varphi}) \boldsymbol{c_{n+1}} = \nabla \mathcal{F}(\boldsymbol{c_n}, \boldsymbol{\varphi}) \boldsymbol{c_n} - \mathcal{F}(\boldsymbol{c_n}, \boldsymbol{\varphi}) \qquad \text{ for every } \boldsymbol{\varphi} \in \boldsymbol{\mathcal{C}}.
\]
\EndWhile
\If{$\|\boldsymbol{c_{n+1}} - \boldsymbol{c_{n}}\| \le \text{tol}$}
\State The algorithm converged.
\State The solution of the variational problem \eqref{eq:vf} is approximated by $\boldsymbol{c_{n+1}}$.
\EndIf
\end{algorithmic}
\end{algorithm}

\begin{figure}[H]
  \centering
  \includegraphics[width=0.5\textwidth]{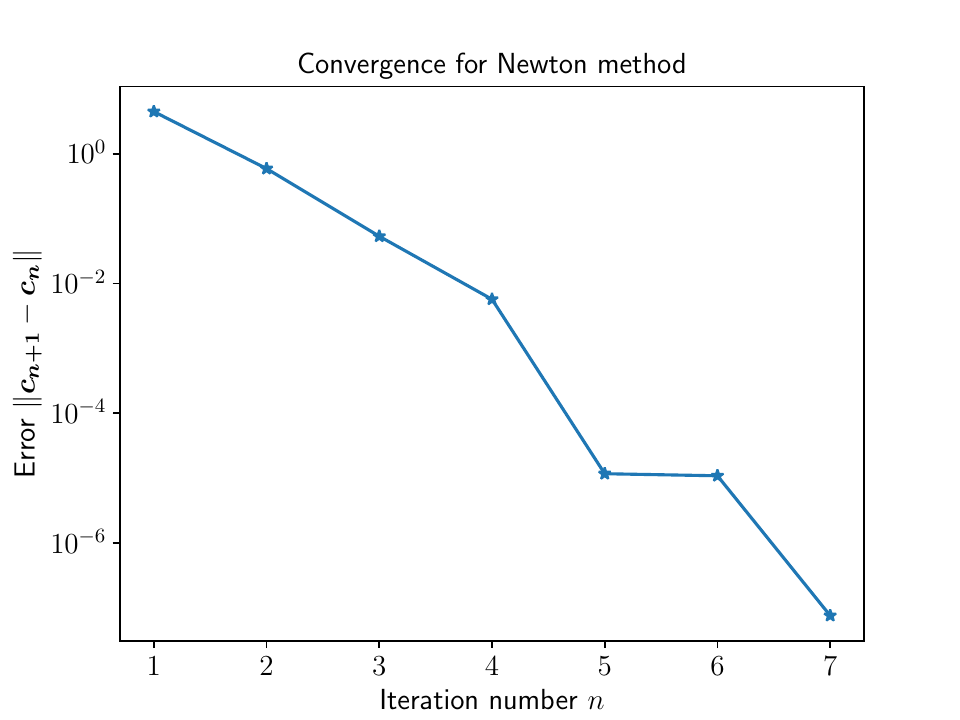}
  \caption{Evolution of the error in Newton algorithm described in Algorithm~\ref{alg:newton} for boundary data given in Table~\ref{tab:bd-p1}.}
  \label{fig:newton-iterates}
\end{figure}

We are now able to precise a few the underlying finite element structure of the fomulation. We consider structured regular triangulation $\mathcal{T}_h$ of $\Omega$ with triangles of diameter $h$ and approach $\mathcal{C}_i$ by its finite dimensional subspace $C_{h}$ formed by the continuous function whose restriction to any triangle of $\mathcal{T}_h$ is a first degree polynomial. All the numerical simulations were done using FreeFEM~\cite{hecht-2012}. In practice, convergence is achieved in a small number of Newton iterations for the parameter regimes considered, and this stationary solver constitutes the forward map used repeatedly within the inverse identification procedure, as described by Figure~\ref{fig:newton-iterates} (please see Table~\ref{tab:bd-p1}).

\section{Inverse Problem Formulation and Numerical Implementation}\label{sec:inverse-pipeline}

We now turn to the numerical resolution of the inverse problem consisting in the identification of effective membrane diffusion coefficients from partial output measurements. The inversion is formulated as a nonlinear least-squares problem, where each evaluation of the cost functional requires the solution of the forward model described in Section~\ref{sec:math-model}.

In order to assess the reliability and stability of the approach, we adopt a progressive validation strategy. We first consider a single-patient configuration, which serves as a baseline validation of the numerical pipeline. We then extend the analysis to a multi-patient setting, allowing us to investigate improvements in parameter identifiability and estimation stability. The robustness of the inversion with respect to noisy observations is subsequently examined, followed by a sensitivity analysis highlighting the respective influence of the diffusion coefficients on the measured outputs. This structured approach enables us to evaluate both the numerical feasibility of the inverse problem and its practical limitations in terms of computational cost and data quality.

\subsection{Single-patient Inverse Problem as a Baseline Validation}
\label{sec:single-patient}

We first consider a single-patient configuration in order to validate the inverse pipeline on a reference dataset. In this setting, inlet concentrations in blood and dialysate, together with outlet blood concentrations, are available from experimental measurements for one patient. These data are subsequently used to generate synthetic datasets for the multi-patient experiments presented in the following sections.

Diffusion within the membrane is modeled through dimensionless coefficients $\alpha_i \in (0,1]$, defined as fractions of the corresponding diffusion coefficients $D_i$ in blood (see \eqref{eq:vf23} in Section~\ref{subsec:forward-stationary}). As discussed in Section~\ref{subsec:transport}, albumin and calcium-albumin complexes do not cross the membrane, and therefore we set $\alpha_2=\alpha_3=0$. Moreover, citrate and calcium-citrate are assumed to exhibit similar diffusion properties, leading to the constraint $\alpha_4=\alpha_5$. As a consequence, only two independent parameters remain to be identified, namely $\alpha_1$ and $\alpha_4$. Let $\overline{\boldsymbol{c}_{\text{input}}} \in \mathbb{R}^{10}$ denote the measured inlet concentrations in blood and dialysate, and let $\overline{\boldsymbol{c}_{\text{output}}} \in \mathbb{R}^5$ be the measured outlet concentrations in blood. For a given vector of diffusion parameters $\boldsymbol{\alpha}=(\alpha_i)_{1\le i\le5}$, the forward model introduced in Section~\ref{subsec:forward-stationary} provides predicted outlet concentrations
\begin{equation*}
	\mathcal{B}_{\boldsymbol{\alpha}}(\overline{\boldsymbol{c}_{\text{input}}}) \in \mathbb{R}^5,
\end{equation*}
obtained by solving the stationary transport problem using FreeFEM. More precisely the output is obtained through the use of the formula \eqref{eq:output-formmula} and the inverse problem is formulated as a nonlinear least-squares minimization. To this end, we introduce the admissible set
\begin{equation*}
	\mathcal{A}=\left\{\boldsymbol{\alpha}\in\mathbb{R}^5 \;|\; \alpha_2=\alpha_3=0,\ \alpha_4=\alpha_5,\ \alpha_i\in(0,1]\ \text{for}\ i\in\{1,4,5\}\right\},
\end{equation*}
and we seek
\begin{equation*}
	\min_{\boldsymbol{\alpha}\in\mathcal{A}} J(\boldsymbol{\alpha}),
\end{equation*}
where the above cost functional is defined by
\begin{equation}\label{eq:J}
	J(\boldsymbol{\alpha})=\sum_{i=1}^5\frac{\big|\mathcal{B}_{\boldsymbol{\alpha},i}(\overline{\boldsymbol{c}_{\text{input}}}) - \overline{c_{i,\text{output}}}\big|^2}{\big|\overline{c_{i,\text{output}}}\big|^2}.
\end{equation}
Due to the previously mentioned constraints, the optimization reduces to two parameters. Introducing $\boldsymbol{\beta}=(\beta_1,\beta_2)$ and setting $\boldsymbol{\alpha}=(\beta_1,0,0,\beta_2,\beta_2)$, the problem becomes the minimization over $[0,1]^2$ of the reduced functional
\begin{equation*}
	\mathcal{J}(\boldsymbol{\beta}) = J(\boldsymbol{\alpha}).
\end{equation*}

In this baseline configuration, the inverse problem is solved directly within the FreeFEM environment using a projected gradient descent method with adaptive step size. For each candidate parameter vector $\boldsymbol{\beta}$, the forward problem is solved using the Newton algorithm described in Algorithm~\ref{alg:newton}, and the resulting outlet concentrations are used to evaluate the reduced cost functional $\mathcal{J}$. Finite-difference approximations of the gradient are computed by repeated calls to the forward solver, and bound constraints are enforced through projection onto $[0,1]^2$. More precisely, for a small parameter $h>0$, we approximate
\begin{equation*}
	\nabla \mathcal{J}(\boldsymbol{\beta})\cdot \boldsymbol{i} \approx \frac{\mathcal{J}(\boldsymbol{\beta}+h\boldsymbol{i})-\mathcal{J}(\boldsymbol{\beta})}{h}, \qquad \nabla \mathcal{J}(\boldsymbol{\beta})\cdot \boldsymbol{j} \approx \frac{\mathcal{J}(\boldsymbol{\beta}+h\boldsymbol{j})- \mathcal{J}(\boldsymbol{\beta})}{h},
\end{equation*}
with $\boldsymbol{i}=(1,0)$ and $\boldsymbol{j}=(0,1)$. The resulting pseudo- gradient vector is then assembled as
\begin{equation}\label{eq:grad-calJ}
	\nabla \mathcal{J}(\boldsymbol{\beta})= \big(\nabla \mathcal{J}(\boldsymbol{\beta})\cdot \boldsymbol{i}, \nabla \mathcal{J}(\boldsymbol{\beta})\cdot \boldsymbol{j}\big)^T.
\end{equation}
The complete optimization procedure is summarized in Algorithm~\ref{alg:grad}. Each evaluation of $\mathcal{J}$ requires the solution of the forward problem, so that several calls to the Newton solver described in Algorithm~\ref{alg:newton} are performed at every iteration. This single-patient implementation serves as a proof of concept for the inverse methodology and provides a reference solution for subsequent experiments. In the following sections, this approach is extended and embedded into a more general FreeFEM-Python computational framework, allowing for multi-patient inversion and the use of more robust derivative-free optimization strategies.

\begin{algorithm}[H]
\caption{Gradient descent algorithm with adaptive step for the minimization of $\mathcal{J}$.}\label{alg:grad}
\begin{algorithmic}
  \Require $\mathcal{J}$, $\boldsymbol{\beta_0}$, $tol$, $n_\text{max}$. \Comment{$n_\text{max}$ is the maximal number of iterations}
  \Require $\text{initial\_step}$ \Comment{$\text{initial\_step}$ is the initial descent step}
  \Require tol.
  \Comment{tol is a small value}
  \State $n \gets 0$.
  \State Compute $\mathcal J(\boldsymbol{\beta_n})$.
  \Comment{Use Newton algorithm described in Algorithm~\ref{alg:newton}}
  
  \State Compute $\nabla \mathcal J(\boldsymbol{\beta_n})$.
  \Comment{Two more calls of Newton method}

  \State $\text{s} \gets \text{initial\_step}$

  \While{$\mathcal{J}(\mathcal{P}(\boldsymbol{\beta_n} - \text{s} \cdot \nabla \mathcal{J}(\boldsymbol{\beta_n}))) \ge \mathcal{J}(\boldsymbol{\beta_n})$ and $\|\mathcal{P}(\boldsymbol{\beta_n} - \text{s} \cdot \nabla \mathcal{J}(\boldsymbol{\beta_n})) - \boldsymbol{\beta_n}\| > \text{tol}$} 
    \State $\text{s} \gets \text{s} / 2$.
    \Comment{$\mathcal{P}(\boldsymbol{\beta})$ is the projection of $\boldsymbol{\beta}$ on $[0, 1]^2$}
  \EndWhile
  \While{$\|\mathcal{P}(\boldsymbol{\beta_n} - \text{s} \cdot \nabla \mathcal{J}(\boldsymbol{\beta_n})) - \boldsymbol{\beta_n}\| > \text{tol}$ and $n < n_\text{max}$}
    \State $\boldsymbol{\beta_{n +1}} \gets \mathcal{P}(\boldsymbol{\beta_n} - \text{s} \cdot \nabla \mathcal{J}(\boldsymbol{\beta_n}))$.
    \State $n \gets n + 1$.
    \State Compute $\nabla \mathcal J(\boldsymbol{\beta_n})$.
    \State $\text{s} \gets \text{initial\_step}$.
      \While{$\mathcal{J}(\mathcal{P}(\boldsymbol{\beta_n} - \text{s} \cdot \nabla \mathcal{J}(\boldsymbol{\beta_n}))) \ge \mathcal{J}(\boldsymbol{\beta_n})$ and $\text{s} \cdot \|\mathcal{P}(\boldsymbol{\beta_n} - \text{s} \cdot \nabla \mathcal{J}(\boldsymbol{\beta_n})) - \boldsymbol{\beta_n}\| > \text{tol}$} 
      \State $\text{s} \gets \text{s} / 2$.
    \EndWhile
  \EndWhile
  \State The approximation of the minima of $\mathcal{J}$ is $\boldsymbol{\beta_n}$.
\end{algorithmic}
\end{algorithm}

As a preliminary numerical illustration, we apply the inverse procedure to the data of one patient from the clinical study. Boundary concentrations at the inlet and outlet of the dialyser are reported in Table~\ref{tab:bd-p1}. 

\begin{table}[ht!]
  \centering
  \begin{tabular}{cccccc}
    \rowcolor{gray!20} Boundary value & $c_1$ & $c_3$ & $c_3$ &$c_4$ & $c_5$ \\
    \cellcolor{gray!20} on $\Gamma_{\ell, b}$ & 0.11 & 3.71602 & 0.0577928 & 5.03048 & 1.37152 \\
    \cellcolor{gray!20} on $\Gamma_{r, d}$ & 1.25 & 0 & 0 & 0 & 0 \\
    \cellcolor{gray!20} on $\Gamma_{r, b}$ & 0.96 & 3.577187 & 0.48553 & 0.144108 & 0.342892
 \\
  \end{tabular}
  \caption{Concentrations (in $\text{mol}\cdot \text{m}^{-3}$) of the five chemical species in blood (boundary $\Gamma_{\ell, b}$) and dialysate (boundary $\Gamma_{r,d}$) at the entrance of the dialyser and in blood at the exit of dialyser (boundary $\Gamma_{r,b}$) for the patient number 1.}
  \label{tab:bd-p1}
\end{table}

Starting from an initial guess $\boldsymbol{\beta}_0=(0.2,0.2)$, the projected gradient descent algorithm converges toward diffusion coefficients that significantly improve the agreement between simulated and measured outlet concentrations.

\begin{figure}[H]
  \centering
  \begin{tabular}{cc}
    \includegraphics[width=0.4\textwidth]{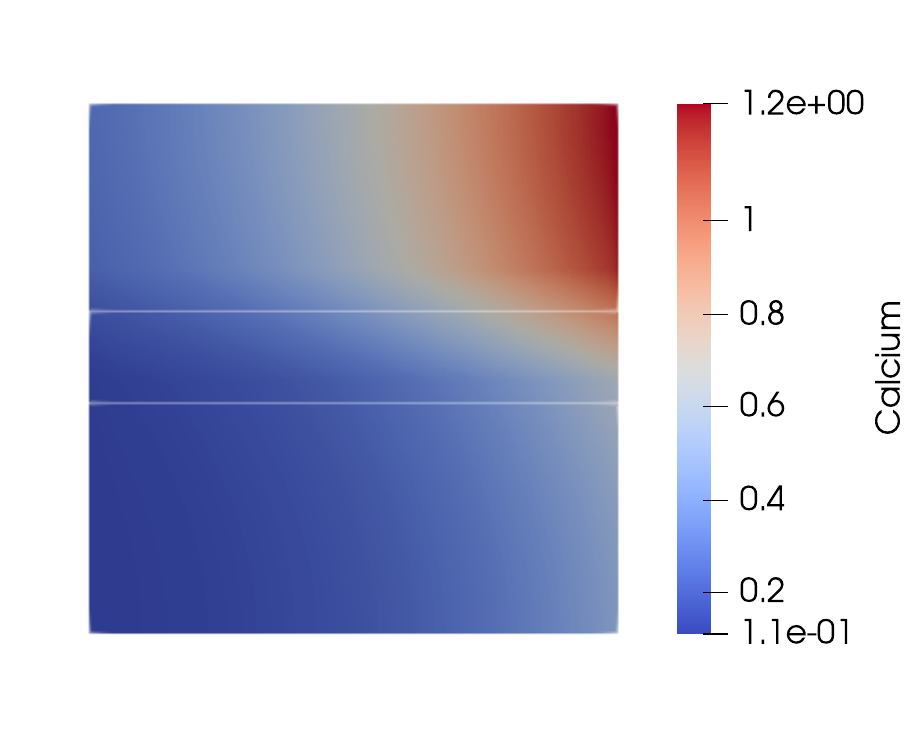}
    & \includegraphics[width=0.4\textwidth]{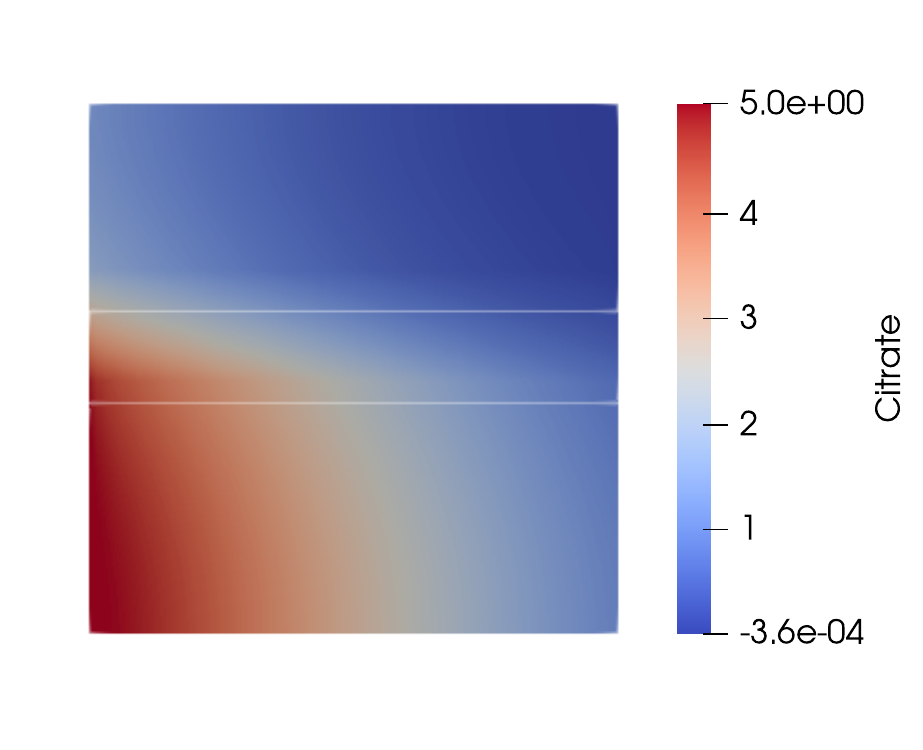}
  \end{tabular}
  \caption{Components $c_1$ (corresponding to Calcium) and $c_4$ (corresponding to Citrate) of the numerical approximation of the solution $\boldsymbol{c}$ computed using Newton's method for the initial data in Table~\ref{tab:bd-p1} and $\boldsymbol{\beta}=(0.2, 0.2)$.} 
  \label{fig:newton-ci}
\end{figure}

\begin{figure}[H]
  \centering
  \begin{tabular}{cc}
    \includegraphics[width=0.4\textwidth]{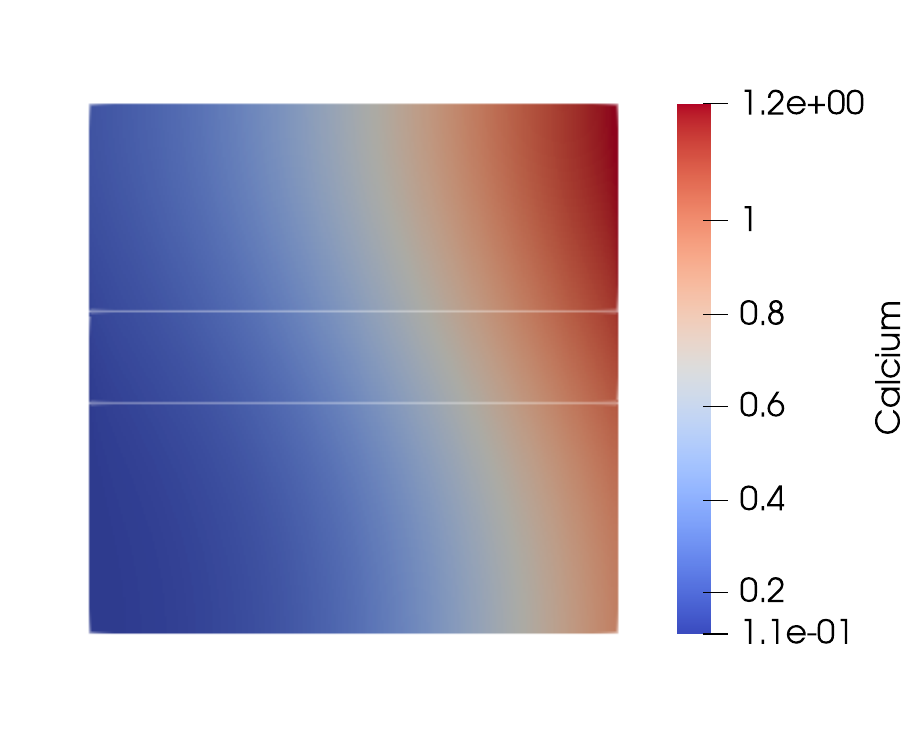}
    & \includegraphics[width=0.4\textwidth]{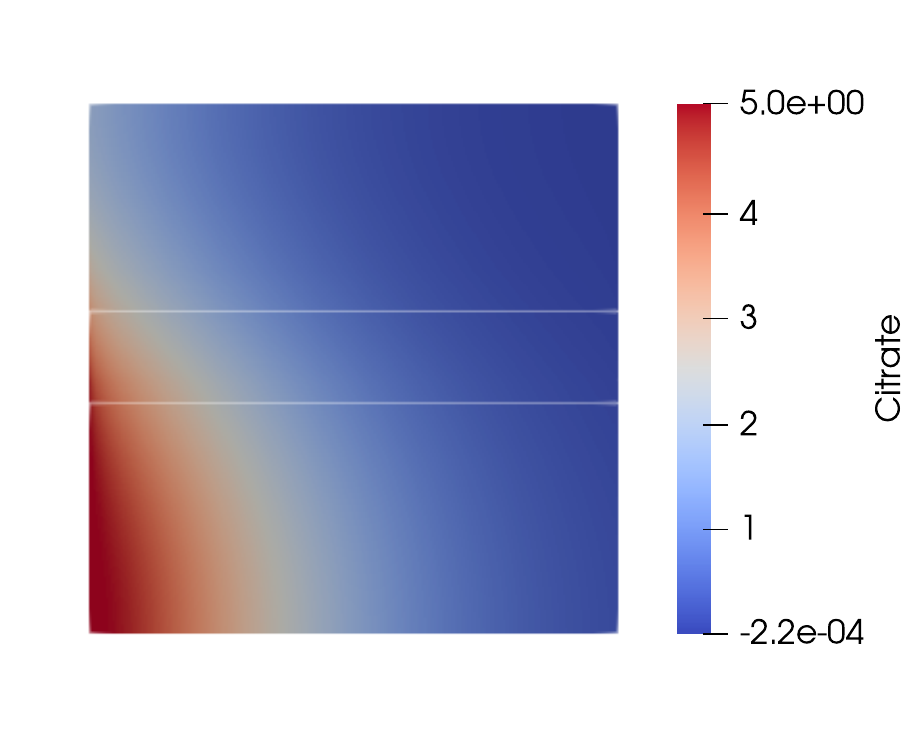}
  \end{tabular}
  \caption{Components $c_1$ and $c_4$ of the numerical approximation of the solution $\boldsymbol{c}$ for $\beta$ minimizing $\mathcal{J}$.} 
  \label{fig:newton-ci-op}
\end{figure}

Figures~\ref{fig:newton-ci} and \ref{fig:newton-ci-op} respectively display the concentrations of Calcium and Citrate for $\beta_0 = (0.2, 0.2)$ and the computed Calcium and Citrate concentrations for the optimized parameters after the convergence of the procedure described above, showing a close match with experimental values at the dialyser outlet. The corresponding value of the cost functional $\mathcal{J}$ is reduced by more than two orders of magnitude compared to the initial guess, indicating that the inverse procedure is able to recover physically meaningful membrane diffusion coefficients in this simplified setting.

While preliminary validation can be performed on a single patient, robust identification of diffusion coefficients requires a multi-patient setting. We therefore construct a synthetic cohort whose statistical structure is consistent with the available clinical data. Beyond inter-patient variability, we further investigate the robustness of the inverse reconstruction with respect to parametric uncertainty by considering two complementary sources of perturbations. First, measurement noise is introduced at the level of the synthetic concentration data, reflecting variability and uncertainty in observable quantities. Second, we assess the sensitivity of the inverse problem with respect to the calcium and citrate diffusion coefficients by perturbing their reference values and analyzing the resulting model outputs.

In this latter setting, a single synthetic patient is selected and multiple forward simulations are performed using perturbed diffusion coefficients of the form $d = d^\ast(1+\varepsilon)$, where $d^\ast$ denotes the nominal value and $\varepsilon$ is a prescribed relative perturbation. The corresponding outputs are then treated as new targets for the inverse procedure. This approach allows us to quantify the stability of the reconstruction with respect to modeling errors in the coefficients, independently of inter-patient variability. Together with the multi-patient experiments based on noisy concentration data, these tests provide a comprehensive numerical assessment of both robustness to measurement noise and sensitivity to parametric perturbations.

\subsection{Multi-patient Inverse Framework and Exact-data Identifiability}\label{subsec:synthetic-identifiability}

Let $\mathcal{D}_{\mathrm{real}}$ denote the table of real patient measurements,
where each row corresponds to a physiological quantity (inlet/outlet pressures,
flows, biochemical concentrations, etc.) and each column corresponds to a patient.
Our goal is to generate a fixed number $N_s$ of synthetic patients while preserving,
for each variable, the empirical mean and standard deviation observed in the real cohort.

For each physiological variable $X$, we proceed independently as follows.
Let $\{x_j\}_{j=1}^{N_r}$ be the available real measurements for that variable,
after removal of missing values. We compute
\begin{equation*}
    \mu_X = \frac{1}{N_r}\sum_{j=1}^{N_r} x_j \quad \text{and}\quad \sigma_X = \sqrt{\frac{1}{N_r-1}\sum_{j=1}^{N_r}(x_j-\mu_X)^2},
\end{equation*}
together with the minimum observed value $m_X = \min_j x_j$. Synthetic samples are then drawn from a normal distribution
\begin{equation*}
    X^{(k)} \sim \mathcal{N}(\mu_X,\sigma_X^2), \qquad k=1,\dots,N_s,
\end{equation*}
and truncated from below by a conservative floor equal to half the minimum observed value,
namely
\begin{equation*}
    X^{(k)} \leftarrow \max\big(X^{(k)},\,0.5\,m_X\big).
\end{equation*}
If $\sigma_X=0$ or undefined, all synthetic values are therefore set equal to $\mu_X$. This simple truncation prevents negative or unphysical values while preserving the first two empirical moments of the real data. All variables are generated independently in such a way, yielding a synthetic patient matrix
\begin{equation*}
    \mathcal{M}_{\mathrm{syn}} \in \mathbb{R}^{M\times N_s},
\end{equation*}
where $M$ is the number of measured fields. A fixed random seed is used to ensure reproducibility. Importantly, only the physiological input data are synthesized at this stage. We emphasize that this synthetic generation ignores cross-variable correlations and is intended only to provide a controlled numerical testbed, not to reproduce full physiological variability.

The dialysis model outputs (five observable targets per patient) are not generated by noise injection. Instead, these are computed consistently by solving the forward problem \eqref{eq:dar145}--\eqref{eq:dar23} for prescribed diffusion coefficients and subsequently used as reference targets in the inverse problem. The procedure implemented in our pipeline is summarized in Algorithm~\ref{alg:generate-synthetic-data}.

\begin{algorithm}\label{alg:generate-synthetic-data}
\caption{Generation of a synthetic multi-patient cohort from real measurements.}\label{alg:synthpatients}
\begin{algorithmic}
  \Require Real patient data table $\mathcal{D}_{\mathrm{real}} \in \mathbb{R}^{M\times N_r}$ (rows = fields, columns = patients).
  \Require Number of synthetic patients $N_s \in \mathbb{N}$.
  \Require Random seed $s$ (optional, for reproducibility).
  \Ensure Synthetic patient data table $\mathcal{M}_{\mathrm{syn}} \in \mathbb{R}^{M\times N_s}$.

  \State Initialize the random number generator with seed $s$.
  \State Initialize $\mathcal{M}_{\mathrm{syn}}$ as an empty $M\times N_s$ table.

  \For{$i = 1,\dots,M$} \Comment{Loop over physiological fields}
    \State Extract the non-missing real values $\{x_j\}_{j=1}^{N_i}$ from row $i$.
    \If{$N_i = 0$}
      \State Set $\mathcal{M}_{\mathrm{syn}}(i,k) \gets \texttt{NaN}$ for all $k=1,\dots,N_s$ and \textbf{continue}.
    \EndIf

    \State Compute $\mu_i \gets \mathrm{mean}(\{x_j\}_{j=1}^{N_i})$.
    \State Compute $\sigma_i \gets \mathrm{std}(\{x_j\}_{j=1}^{N_i})$ (sample standard deviation).
    \State Compute $m_i \gets \min(\{x_j\}_{j=1}^{N_i})$.
    \State Set the lower bound $\ell_i \gets 0.5\,m_i$. \Comment{Conservative floor to avoid unphysical values}

    \If{$\sigma_i$ is undefined or $\sigma_i = 0$}
      \State Set $y_k \gets \mu_i$ for all $k=1,\dots,N_s$. \Comment{Constant row}
    \Else
      \For{$k = 1,\dots,N_s$}
        \State Draw $y_k \sim \mathcal{N}(\mu_i,\sigma_i^2)$.
      \EndFor
    \EndIf

    \For{$k = 1,\dots,N_s$}
      \State Enforce the floor: $y_k \gets \max(y_k,\ell_i)$.
      \State Set $\mathcal{M}_{\mathrm{syn}}(i,k) \gets y_k$.
    \EndFor
  \EndFor

  \State \Return $\mathcal{M}_{\mathrm{syn}}$.
\end{algorithmic}
\end{algorithm}

Also, in addition to sampling physiological input variables, a calibration step is performed to construct a consistent set of synthetic reference data. More precisely, for each synthetic patient, auxiliary boundary quantities (in particular pressure profiles or hydraulic parameters) are determined by solving an inverse problem within the forward pipeline. This calibration ensures that the forward model reproduces prescribed reference concentration levels at the dialyzer inlet and outlet, thereby enforcing physical consistency between inputs and outputs. This preliminary inverse step is not the focus of the present work and is only used to generate coherent synthetic datasets. The calibrated pressure data are subsequently treated as fixed patient-specific inputs, while the diffusion coefficients remain unknowns in the multi-patient inverse problem considered below.

Thus, such a synthetic cohort is subsequently used to generate patient-specific forward
solutions of the dialysis model, from which the observable outputs are extracted. These outputs constitute the multi-patient targets employed in the inverse problem for diffusion coefficient identification.

We now turn to the identification of the diffusion coefficients in a multi-patient setting,
using the synthetic cohort constructed above. While single-patient experiments provide a first validation of the modeling and inversion pipeline, aggregating information across several patients is needed to improve robustness with respect to noise and to mitigate patient-specific variability.

Let $\mathcal{P}=\{1,\dots,N_p\}$ denote the set of synthetic patients, namely the columns of $\mathcal{M}_{\mathrm{syn}}$ described above. For each patient $p\in\mathcal{P}$ and for a given vector of diffusion coefficients $\boldsymbol{\beta}=(d_{\mathrm{Ca}},d_{\mathrm{Ci}}) \in (0,+\infty)^2$, we denote by $F_p(\boldsymbol{\beta})\in\mathbb{R}^m$ the vector of observable quantities produced by the forward solver, where $m=5$ in the present study.
The mapping $F_p$ is obtained by solving the convection-reaction-diffusion model with
patient-specific inputs (geometry, calibrated pressure data, and synthetic physiological
parameters), and extracting the corresponding outlet quantities at the end of the simulation. From now on, we should denote by $y_p^{\mathrm{obs}}\in\mathbb{R}^m$ the reference target vector associated to patient $p$, generated as described in the previous section.

The diffusion coefficients are estimated by minimizing a global least-squares misfit that
aggregates residuals over all patients:
\begin{equation}\label{eq:Jmulti}
    \mathcal{J}(\boldsymbol{\beta}) = \sum_{p\in\mathcal{P}}\|W\cdot\big(F_p(\boldsymbol{\beta})-y_p^{\mathrm{obs}}\big)\|_2^2
    + \lambda\mathcal{R}(\boldsymbol{\beta}).
\end{equation}
Here $W \in \mathbb{R}^{m \times m}$ ($m = 5$) is a diagonal scaling operator used to normalize the different components of the output vector, so that no single observable dominates the objective due to its physical units or magnitude. Usually, the diagonal entries of $W$ are chosen as the inverse of empirical scales (for example mean absolute values or standard deviations computed over the cohort). The optional regularization term $\mathcal{R}$ enforces weak a priori constraints on the parameter range, and $\lambda \geq 0$ denotes its weight. In the synthetic experiments reported below, $\lambda$ is taken small or equals to zero whenever the inverse problem is sufficiently well-conditioned.

The evaluation of $\mathcal{J}$ requires solving the forward system \eqref{eq:dar145}--\eqref{eq:dar23} for each patient and is therefore computationally expensive. Moreover, the forward solver is treated as a ``black box'', so that analytic gradients are unavailable. For these reasons, we adopt a derivative-free direct-search strategy to minimize the functional $\mathcal{J}$ given in \eqref{eq:Jmulti}. Since there is only two parameters to estimate, classical methods such as Powell algorithm are well suited and typically converge in a moderate number of objective evaluations. Powell method proceeds by iteratively performing one-dimensional line searches along a set of directions in parameter space. Starting from an initial guess, the objective function is minimized successively along each direction, and at the end of each iteration, a new direction is constructed from the net displacement between successive iterates and replaces the least effective previous one. This mechanism allows the method to progressively build approximate curvature information without computing derivatives. In practice, Powell algorithm is expected to converge in a moderate number of objective evaluations and provides stable reconstructions in both noise-free and noisy experiments, we refer the interested reader to the book of A.Conn, K.Scheinberg \& L.Vicente~\cite{conn-scheinberg-vicente-2009} for more details concrning Powell-type and more generally derivative-free methods.

Also, to enforce positivity and improve numerical conditioning, the parameters are reparametrized as $\boldsymbol{\beta}=(e^z_1,e^z_2)$, and the optimization is performed with respect to $z\in\mathbb{R}^2$. Admissible bounds on $\boldsymbol{\beta}$ are imposed through a penalty term in~\eqref{eq:Jmulti}. Forward failures or nonphysical parameter values are handled robustly by returning a large objective value, allowing the optimizer to automatically avoid such regions. In the single-patient setting considered previously, a pseudo-gradient descent with adaptive step size has been efficiently implemented, as finite-difference approximations of the gradient remain affordable. In contrast, in the multi-patient framework each evaluation of $\mathcal{J}$ already involves $N_s$ forward solves, and numerical gradients would multiply this cost. Derivative-free optimization therefore provides a more favorable compromise between robustness and computational efficiency in the present low-dimensional parameter setting.

Also, as aforementioned a major computational issue lies in the evaluation of $\mathcal{J}(\boldsymbol{\beta})$, which requires computing the patient-specific forward maps $(F_p(\boldsymbol{\beta}))_{p\in\mathcal{P}}$. Since these forward problems are mutually independent, we parallelize their computation on a
single workstation by launching multiple patient simulations concurrently using a local multi-process strategy. This reduces the wall-clock time of each objective evaluation almost proportionally to the number of available CPU cores, up to input and output limitations. In addition, optional heavy outputs (such as visualization files) are disabled during the inverse procedure, as they are not required for objective evaluation and significantly impact runtime. The above-mentioned reparametrization $\boldsymbol{\beta}=(e^{z_1},e^{z_2})$ is
used throughout the optimization, and admissible parameter ranges are enforced
through a soft quadratic penalty added to the objective functional. Forward
solver failures as well as nonphysical parameter values are handled by returning a large objective value. 

While the synthetic cohort provides physiological input parameters, additional
boundary quantities required by the forward  system \eqref{eq:dar145}--\eqref{eq:dar23}, notably the inlet and outlet pressures on both blood and dialysate sides, cannot be prescribed independently, as they must satisfy \eqref{eq:fluid-flow}. For this reason, we recall once again that a preliminary calibration step is carried out for each synthetic patient. In this auxiliary inverse problem, the pressure values at the inlet and outlet boundaries are treated as unknowns and determined so that the forward system reproduces prescribed reference concentration levels, this step ensuring consistency between flow conditions and transport dynamics.

Concretely, for each synthetic patient, the physiological parameters sampled
from the real cohort are first fixed. An auxiliary inverse hydraulic problem is
then solved to identify inlet and outlet pressures consistent with the model and
prescribed reference concentration levels, as aforementioned. These calibrated pressures are subsequently frozen as patient-specific inputs, and the full convection–reaction–diffusion model is finally solved with prescribed diffusion coefficients $(d_{\mathrm{Ca}},d_{\mathrm{Ci}})$ to generate the observable
outputs, so that the resulting outlet concentrations constitute the exact synthetic reference targets $y_p^{\mathrm{obs}}$ used in the inverse study. This two-stage procedure, namely an hydraulic calibration followed by a forward transport simulation, guarantees physical coherence between sampled physiological inputs, boundary conditions, and model outputs. Importantly, the synthetic observables are not only generated by direct sampling or artificial noise injection, but are obtained by solving the complete system \eqref{eq:dar145}--\eqref{eq:dar23} after patient-specific calibration. 

\begin{figure}
	\begin{tabular}{cc}
	\includegraphics[width=0.45\textwidth]{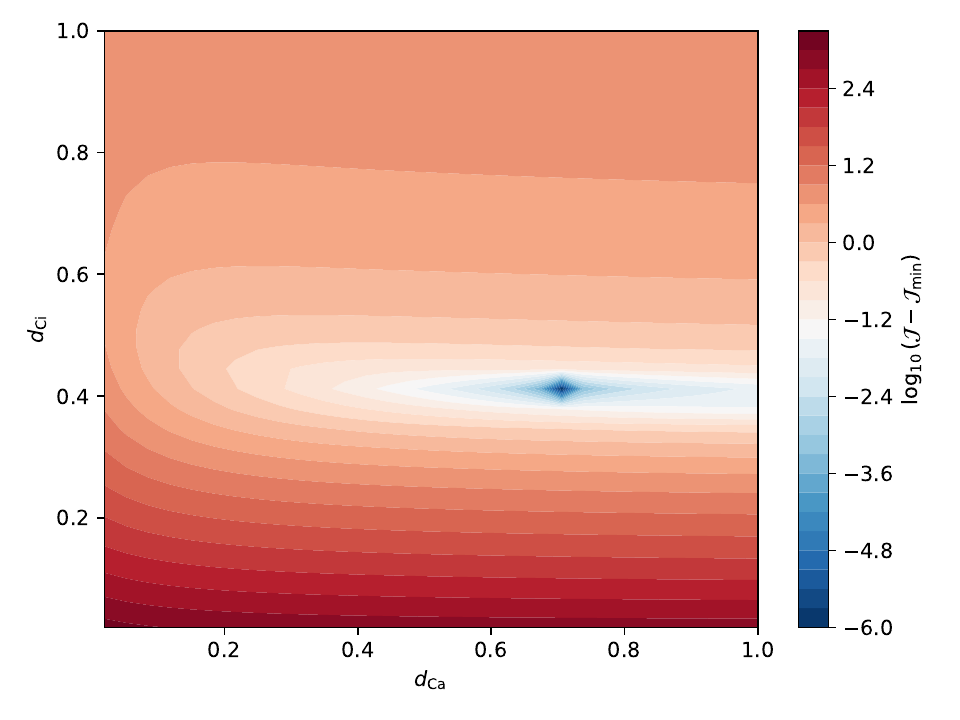} & \includegraphics[width=0.55\textwidth]{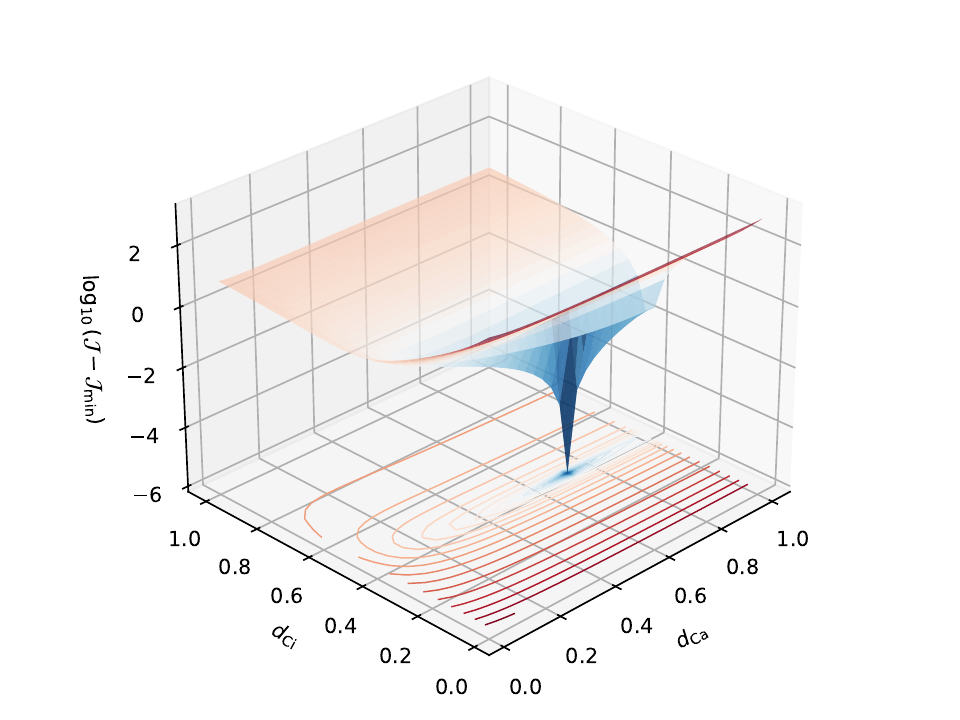}\\
	(a) & (b)
	\end{tabular}
	\caption{The geometry of the functional $\mathcal{J}$ as a function of $\beta = (d_{\mathrm{Ca}},d_{\mathrm{Ci}})$ in logarithmic scale.}
	\label{fig:identifiability-geometry}
\end{figure}

Following such a procedure, we generate a synthetic cohort of $N_s=40$ patients as described above. For each synthetic patient, patient-specific boundary conditions are first calibrated, and the five observable outputs are then generated by solving the forward system \eqref{eq:dar145}--\eqref{eq:dar23} with a fixed ground-truth diffusion vector
\begin{equation*}
	\boldsymbol{\beta}^\star=(d_{\mathrm{Ca}}^\star,d_{\mathrm{Ci}}^\star)=(0.8,0.4).
\end{equation*}
These model-consistent outputs are stored as exact reference targets $\{y_p^{\mathrm{obs}}\}_{p\in\mathcal{P}}$. Before launching the derivative-free optimization, we perform a coarse exploratory grid search of the reduced objective functional $\mathcal J(\beta)$ over the domain $[0.02,1]\times[0.02,1]$, where $\beta=(d_{\mathrm{Ca}},d_{\mathrm{Ci}})$. The purpose of this preliminary step is not to obtain an accurate minimizer, but rather to visualize the global geometry of $\mathcal J$ and to localize a priori regions of potential local or global minima. Figure~\ref{fig:identifiability-geometry} displays $\mathcal J$ in logarithmic scale on a $31 \times 31$ uniform grid. A pronounced valley structure clearly emerges, centered around $d_{\mathrm{Ci}}\approx 0.4$, while the functional exhibits a much flatter dependence with respect to $d_{\mathrm{Ca}}$. In particular, the deepest part of the valley is observed for $0.6 \lesssim d_{\mathrm{Ca}} \lesssim0.9$, indicating a strong sensitivity of the misfit to citrate diffusion, contrasted with a weaker and more correlated influence of calcium one. This anisotropic geometry suggests that $d_{\mathrm{Ci}}$ is more tightly identifiable than $d_{\mathrm{Ca}}$, at least in the present illustrating example, the latter lying along a relatively flat direction of the objective landscape. Such a structure is typical of partially correlated parameters and explains why different optimization runs may converge to slightly different values of $d_{\mathrm{Ca}}$ while remaining close in $d_{\mathrm{Ci}}$. Beyond providing physical insight into parameter sensitivity, this grid search also serves as a diagnostic tool to guide the initialization of Powell method and to mitigate the risk of convergence toward spurious local minima.

To assess numerical identifiability in a multi-patient setting, we then solve the inverse problem~\eqref{eq:Jmulti} using a subset of four patients (here $p\in\{s12,s13,s14,s15\}$). The objective function aggregates the least-squares residuals over these patients, and the resulting two-parameter optimization problem is minimized with the Powell method described above. Convergence is monitored through the sequence of iterates $(d_{\mathrm{Ca}}^{(k)},d_{\mathrm{Ci}}^{(k)})$ and the corresponding objective
values $\mathcal{J}^{(k)}$.

\begin{figure}[H]
	\centering
	\includegraphics[width=0.5\textwidth]{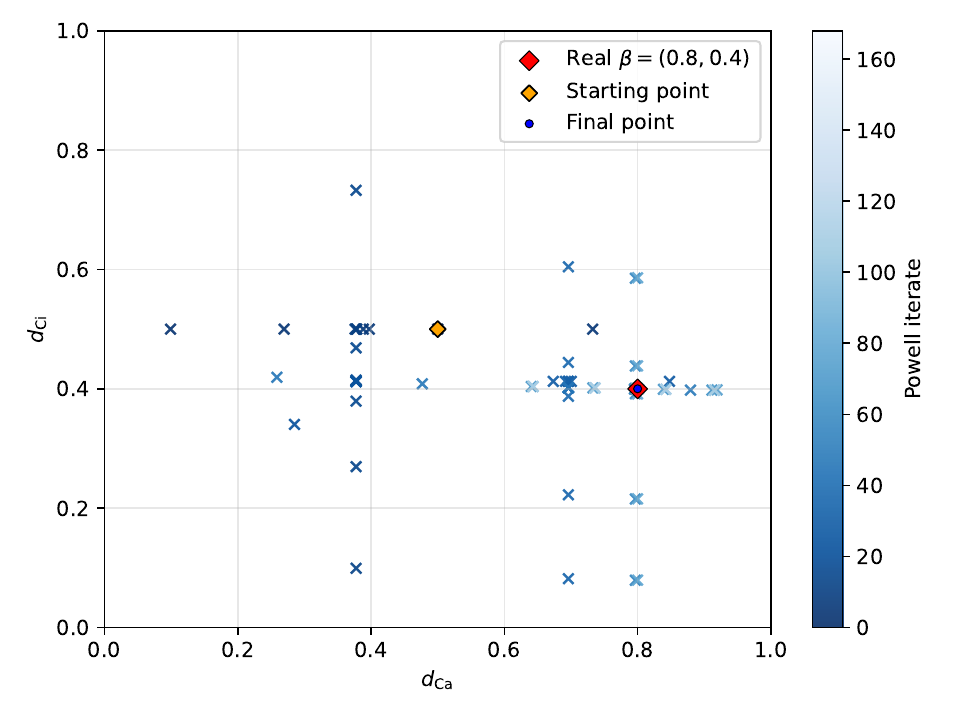}
	\caption{Convergence of the Powell method in the phase plane for the prescribed data $\boldsymbol{\beta} = (0.8,0.4)$ without noise}
	\label{fig:identifiability-powell-2d}
\end{figure}

Figure~\ref{fig:identifiability-powell-2d} illustrates the behavior of Powell algorithm in parameter space for the exact-data multi-patient experiment. The successive iterates $(d_{\mathrm{Ca}}^{(k)},d_{\mathrm{Ci}}^{(k)})$ are displayed in the $(d_{\mathrm{Ca}},d_{\mathrm{Ci}})$ plane, together with the ground-truth value $\boldsymbol{\beta}^\star=(0.8,0.4)$. Starting from an initial guess, the optimization trajectory rapidly enters a neighborhood of the true parameters and subsequently converges toward $\boldsymbol{\beta}^\star$.

\begin{figure}[H]
	\centering
	\includegraphics[width=0.5\textwidth]{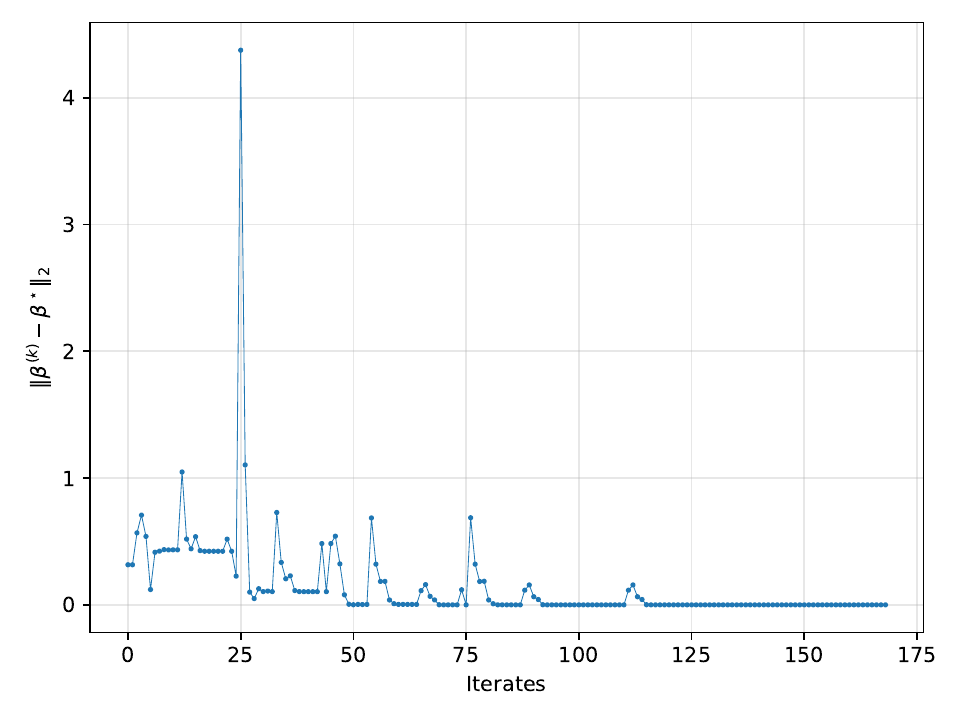}
	\caption{Error of convergence of the Powell method for the prescribed data $\boldsymbol{\beta} = (0.8,0.4)$ without noise}
	\label{fig:identifiability-powell-error}
\end{figure}

Figure~\ref{fig:identifiability-powell-error} shows the evolution of the Euclidean distance $\|\boldsymbol{\beta}^{(k)}-\boldsymbol{\beta}^\star\|_2$ between the current estimate and the ground-truth diffusion coefficients along the Powell iterations. A rapid decrease of the parameter error is observed during the early stages of the optimization, followed by a slower refinement phase as the iterates approach the exact solution. The error eventually reaches numerical precision, confirming convergence of the reconstructed parameters toward the true values. 

\begin{figure}[H]
	\centering
	\includegraphics[width=0.5\textwidth]{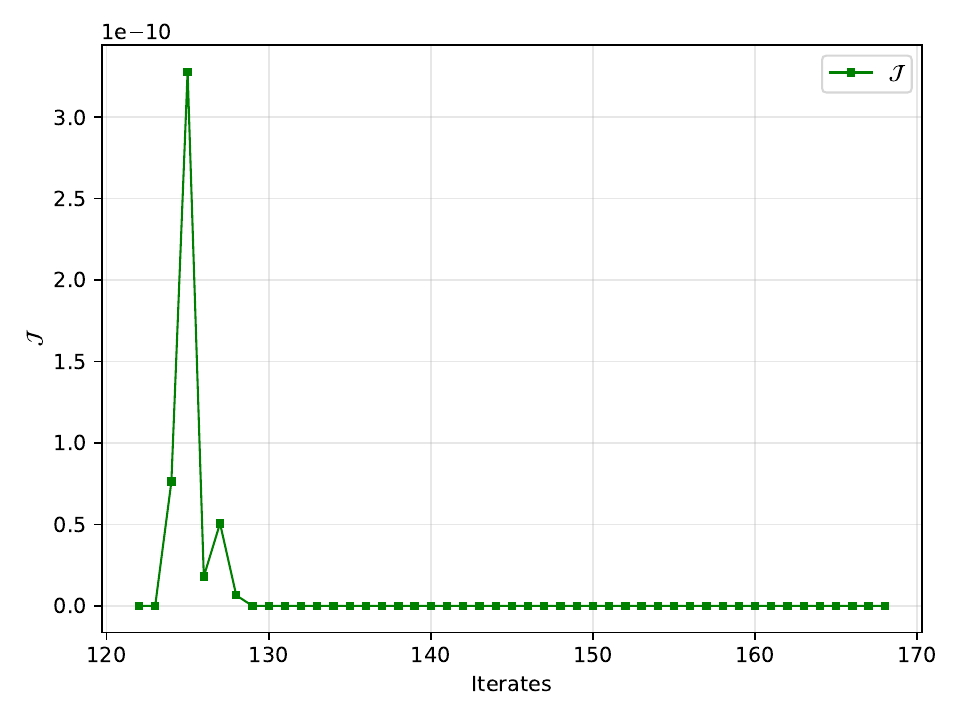}
	\caption{Convergence of the functional for the Powell method with prescribed data $\boldsymbol{\beta} = (0.8,0.4)$ without noise}
	\label{fig:identifiability-powell-J}
\end{figure}

The corresponding evolution of the last iterations of the objective functional $\mathcal{J}$ is shown in Figure~\ref{fig:identifiability-powell-J}. A rapid decrease of the misfit is observed during the first iterations, followed by a slower refinement phase, until numerical convergence is reached, as illustrated in the Figure. Both the parameter trajectory and the objective values indicate
stable convergence toward the exact solution.

These results highlight that, in the absence of measurement noise, the multi-patient inverse problem admits a well-resolved minimum at the true diffusion coefficients. In other words, aggregating information from several patients restores numerical identifiability of $(d_{\mathrm{Ca}},d_{\mathrm{Ci}})$ within the proposed framework. This experiment provides a baseline validation of the multi-patient pipeline before investigating robustness with respect to measurement noise and parametric uncertainty.

\subsection{Robustness with respect to Measurement Noise}\label{subsec:synthetic-robustness}

After establishing numerical identifiability in the exact-data setting, we now
investigate the robustness of the multi-patient inverse problem with respect to
measurement noise on the observable outputs. Starting from the synthetic cohort constructed in Section~\ref{subsec:synthetic-identifiability}, we generate three
additional datasets by perturbing the reference outlet concentrations with
artificial noise. Only the observable targets are modified, while all physiological inputs, calibrated boundary conditions, and model parameters are
kept unchanged. This allows us to isolate the effect of measurement uncertainty
on the inverse reconstruction.

Let $y_p^{\mathrm{obs}}\in\mathbb{R}^m$ denote the exact target vector associated with patient $p$, where $m=5$ in the present study. For each component, noisy observations are generated according to a multiplicative perturbation model
\begin{equation*}
	y_{p,i}^{\mathrm{noisy}} = y_{p,i}^{\mathrm{obs}}\,(1+\varepsilon_{p,i}),
\end{equation*}
where the random variables $\varepsilon_{p,i}$ are drawn independently from a
centered normal distribution with prescribed standard deviation $\sigma$. To
avoid unrealistically large perturbations, the samples are clipped to a bounded interval proportional to $\sigma$. Also, we generate the synthetic cohorts by considering three diffrent noise levels, corresponding to relative standard deviations of $1\%$, $3\%$, and $5\%$. Namely, for each such noise level a complete noisy synthetic cohort is generated, yielding three independent inverse test cases. A fixed random seed is used to ensure reproducibility. The inverse problem is then solved using the same multi-patient objective functional and optimization strategy as in the exact-data case.

\begin{figure}[H]
  \centering
  \includegraphics[width=0.7\textwidth]{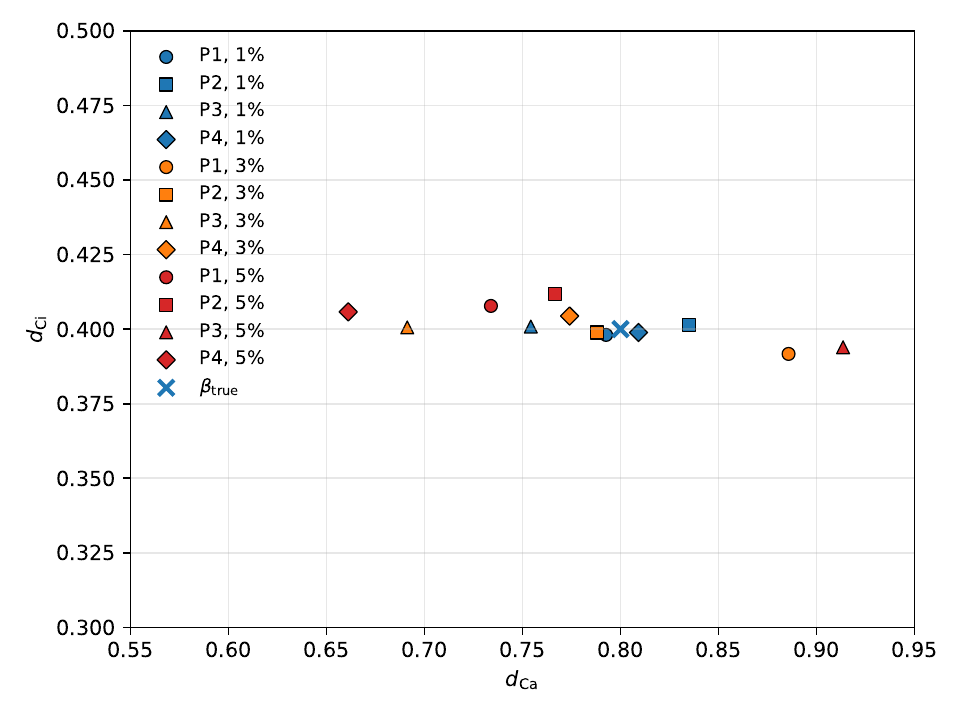}
  \caption{Estimated coefficients $\beta^\star$ for four disjoint sub-cohorts of five synthetic patients, for noise levels $1\%$, $3\%$, and $5\%$ on the targets.}
  \label{fig:synthetic-noise-subcohorts-scatter}
\end{figure}

This procedure mimics measurement uncertainty at the level of clinically
observable quantities while preserving full physical consistency of the forward
model. It provides a controlled framework to assess the stability of the
estimated diffusion coefficients with respect to increasing noise amplitude. For each noise level, we solve the same multi-patient inverse problem as in the exact-data case, using Powell derivative-free method and identical bounds for $(d_{\mathrm{Ca}},d_{\mathrm{Ci}})$. In order to probe the stability of the reconstruction with respect to both the noise amplitude and the cohort composition, we consider four disjoint sub-cohorts of five patients, namely $(s1\text{-}s5)$, $(s6\text{-}s10)$, $(s11\text{-}s15)$, and $(s16\text{-}s20)$, that we denoted by $P1$, $P2$, $P3$ and $P4$ respectively. For each sub-cohort and each noise level, we report the estimated coefficient vector $\beta^\star=(d_{\mathrm{Ca}}^\star,d_{\mathrm{Ci}}^\star)$. Figure~\ref{fig:synthetic-noise-subcohorts-scatter} summarizes the results in parameter space. Overall, the reconstructions remain close to the ground-truth coefficient $\beta_{\mathrm{true}}$ even when the observable targets are perturbed. As expected from the objective landscape discussed in Section~\ref{subsec:synthetic-identifiability}, the estimates exhibit a stronger dispersion along the $d_{\mathrm{Ca}}$-direction than along the $d_{\mathrm{Ci}}$-direction, consistently with a weaker identifiability of $d_{\mathrm{Ca}}$.

\begin{figure}[H]
  \centering
  \includegraphics[width=0.7\textwidth]{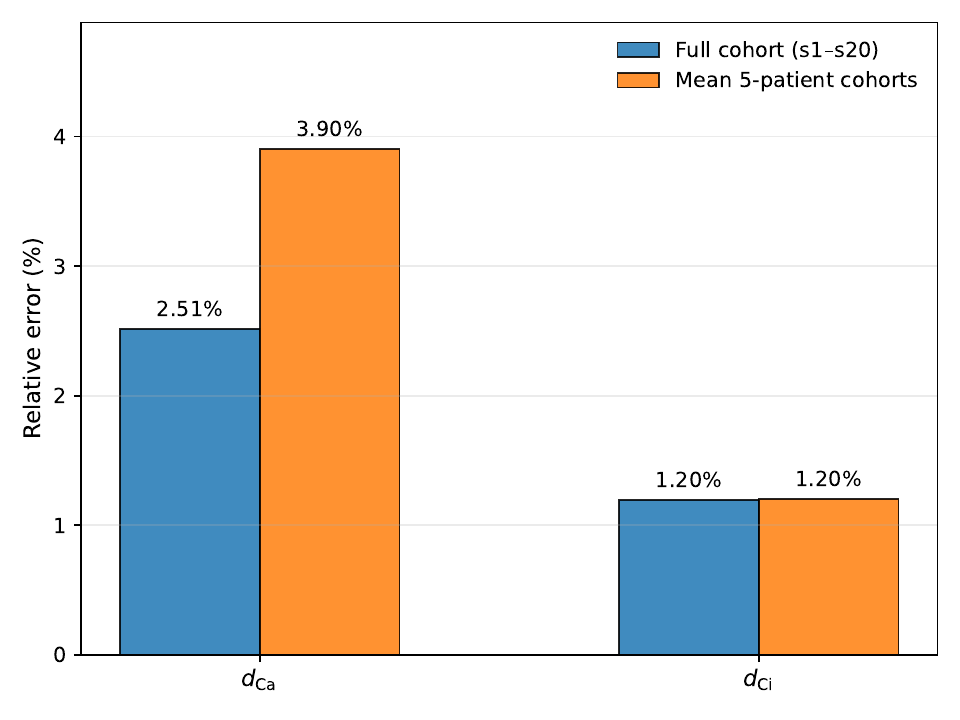}
  \caption{Relative error in percents of the reconstructed coefficients for the full cohort $(s1\text{-}s20)$ compared to the error of means for sub-cohorts with five synthetic patients in the case of $5\%$ noise on the targets.}
  \label{fig:synthetic-noise-s1s20-errors}
\end{figure}

To further assess robustness at a larger scale, we also solve the inverse problem for the full cohort $(s1\text{-}s20)$ using the $5\%$ noisy targets. The recovered coefficient $\beta^\star_{20}$ remains close to $\beta_{\mathrm{true}}$, and Figure~\ref{fig:synthetic-noise-s1s20-errors} reports the resulting relative errors on each component of $\beta$, compared to the mean of errors for five-patients synthetic cohorts over the same patients, namely the means of errors of $P1$, $P2$, $P3$ and $P4$. This confirms that, in this controlled setting, measurement noise at the outlet level does not induce an excessive amplification in the reconstructed effective diffusion parameters. Moreover, it ihghlights that more patients are taken into account the functional; lower is the noise over the diffusion coefficient of Calcium, due to its weaker identifiability. It confirms that such a multi-patient setting is particularly suited for such estimation of the diffusion coefficients. 

While these experiments already provide meaningful evidence of robustness on cohorts of up to $N = 20$ patients, they remain computationally expensive because each functional evaluation requires solving the forward model for all patients in the cohort. This cost prevents a systematic exploration of larger cohorts, as the full synthetic cohort of $N_s = 40$ synthetic patients, and more extensive uncertainty quantification through repeated noisy realizations. This motivates the development of computationally efficient surrogate models, which we discuss in Section~\ref{sec:conclusions}. Finally, the above robustness results concern uncertainty in the observed targets. In the next subsection, we investigate the complementary question of parameter sensitivity, namely how perturbations of $(d_{\mathrm{Ca}},d_{\mathrm{Ci}})$ propagate through the forward model and affect the outlet observables.

\subsection{Sensitivity with respect to Membrane Diffusion Coefficients}\label{subsec:synthetic-sensitivity}

We finally investigate the intrinsic sensitivity of the forward model with respect to perturbations of the membrane diffusion coefficients. This analysis aims at quantifying how uncertainties in $(d_{\mathrm{Ca}}, d_{\mathrm{Ci}})$ propagate through the nonlinear convection-reaction-diffusion system and affect the observable outlet concentrations. Starting from the nominal diffusion vector
\begin{equation*}
	\beta^\star = (d_{\mathrm{Ca}}^\star,d_{\mathrm{Ci}}^\star)=(0.8,0.4),
\end{equation*}
as used above, we generate perturbed coefficients for each of the $N_s=40$ synthetic patients according to
\begin{equation*}
	d = d^\star(1+\varepsilon),
\end{equation*}
where $\varepsilon$ is drawn independently from a centered normal distribution with relative standard deviations of $1\%$, $3\%$, and $5$, following the same calculations as previously done to generate such a noise. Then for each noise level, the full forward problem \eqref{eq:dar145}--\eqref{eq:dar23} is solved with these modified coefficients, while all physiological inputs and calibrated hydraulic boundary conditions are kept fixed. The resulting outlet concentrations are written as new synthetic targets. Importantly, noise is not injected directly at the level of observables. Instead, the perturbations are applied solely to $(d_{\mathrm{Ca}},d_{\mathrm{Ci}})$ and propagated through the complete forward solver, ensuring physical consistency of the generated data, that are represented in Figure~\ref{fig:noisy-coeff} below.

\begin{figure}[H]
	\centering
	\begin{tabular}{cc}
		\includegraphics[width=0.5\textwidth]{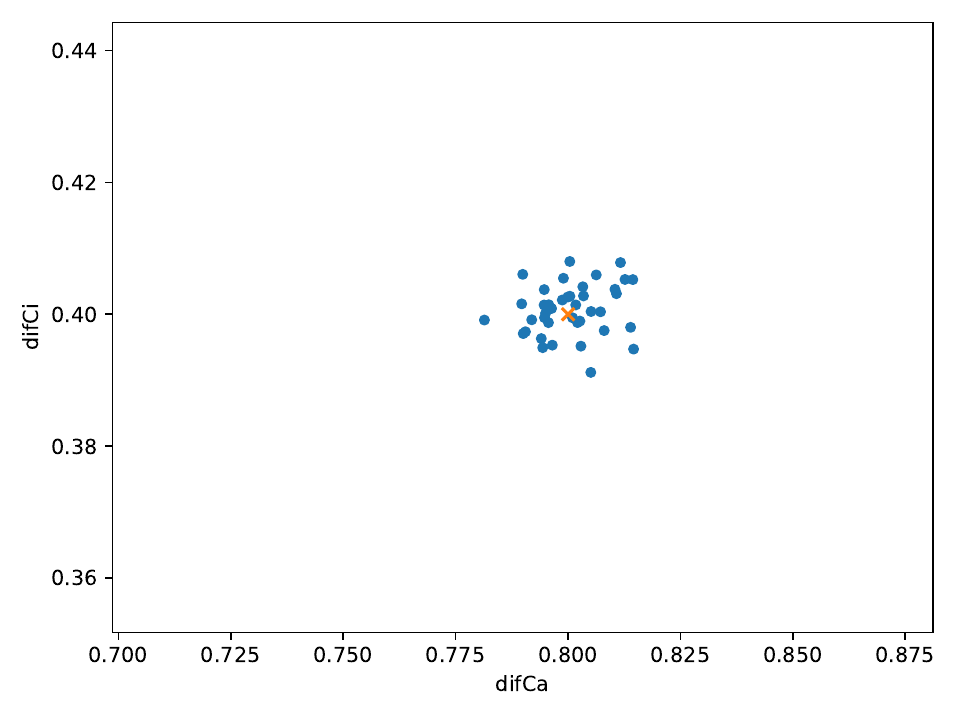}& \includegraphics[width=0.5\textwidth]{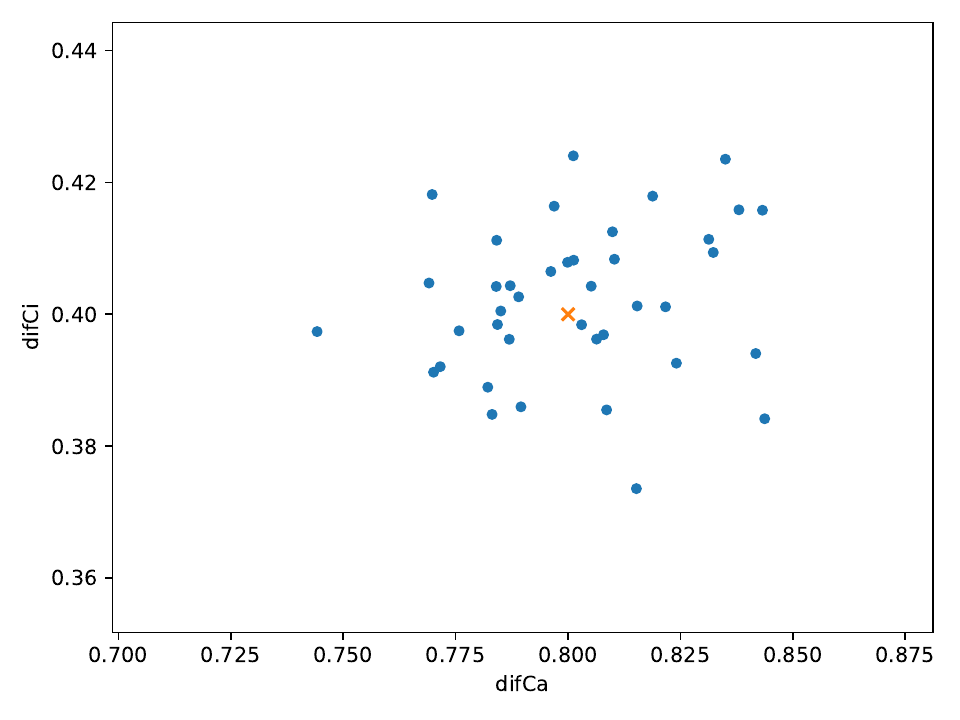}\\
		(a)&(b)
	\end{tabular}\\
	\begin{tabular}{c}
		\includegraphics[width=0.5\textwidth]{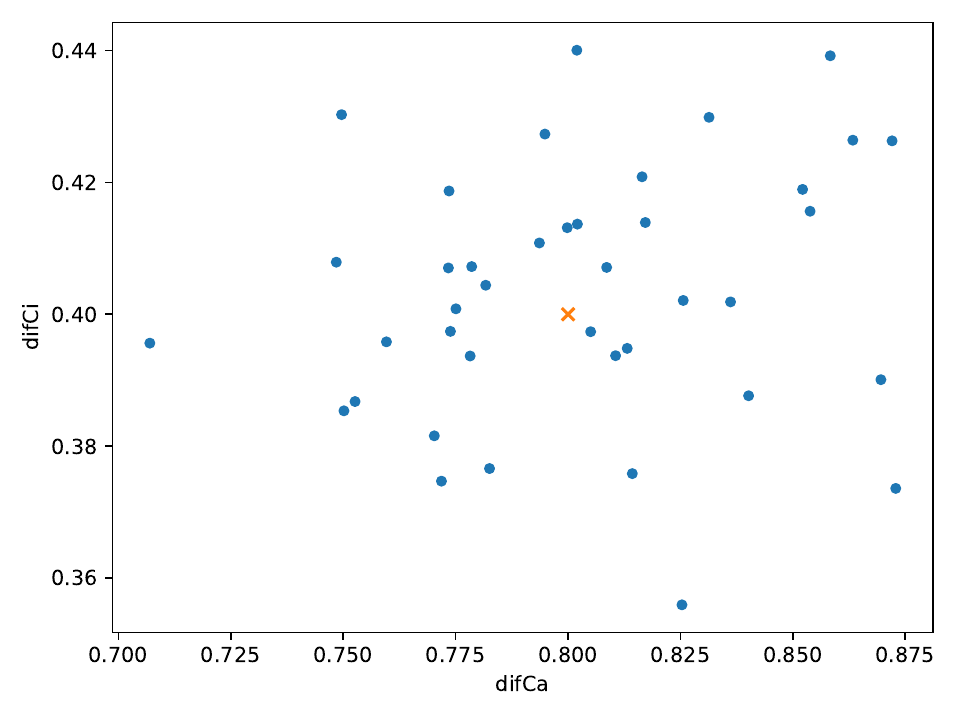}\\
		(c)
	\end{tabular}
	\caption{Noisy diffusion coefficients for synthetic cohort with a noise of order (a) 1\% (b) 3\% and (c) 5\%.} 
	\label{fig:noisy-coeff}
\end{figure}

\begin{figure}[H]
	\centering
	\begin{tabular}{cc}
		\includegraphics[width=0.5\textwidth]{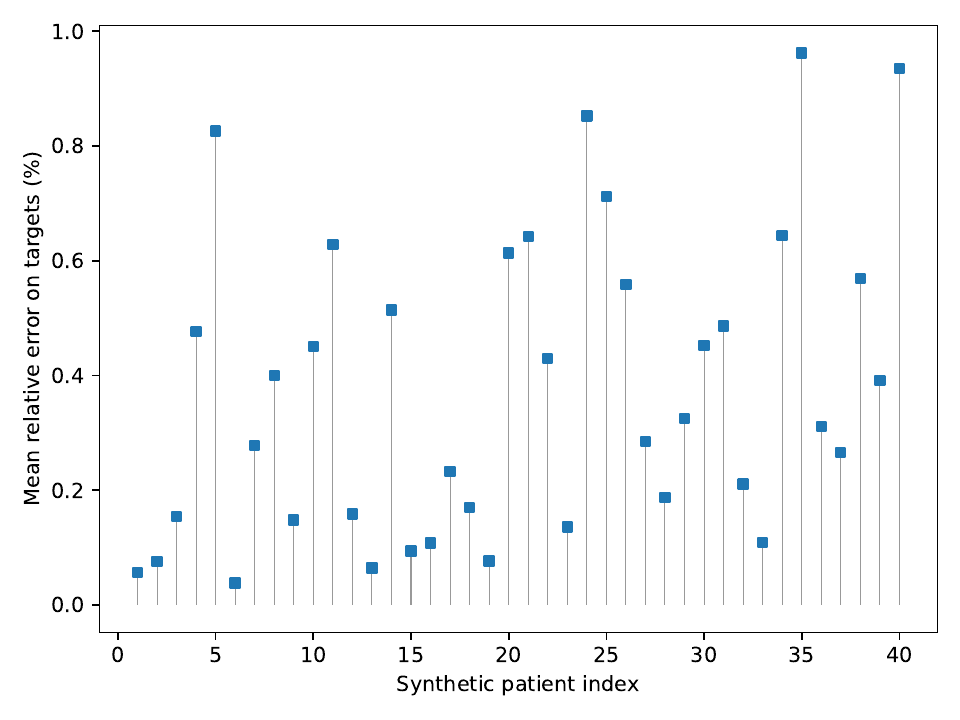}& \includegraphics[width=0.5\textwidth]{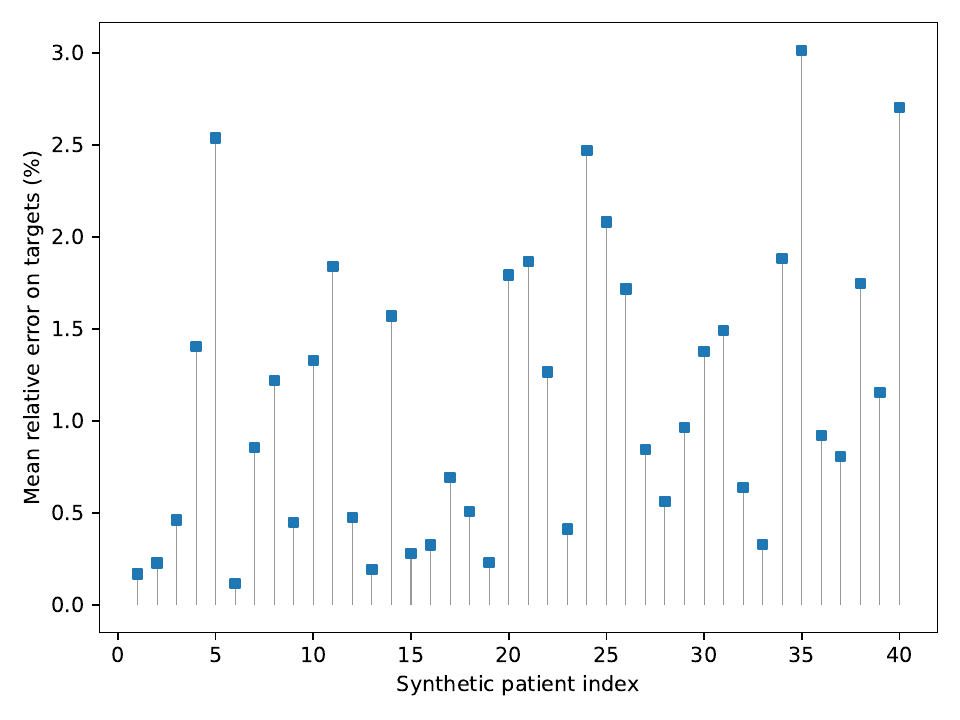}\\
		(a)&(b)
	\end{tabular}\\
	\begin{tabular}{c}
		\includegraphics[width=0.5\textwidth]{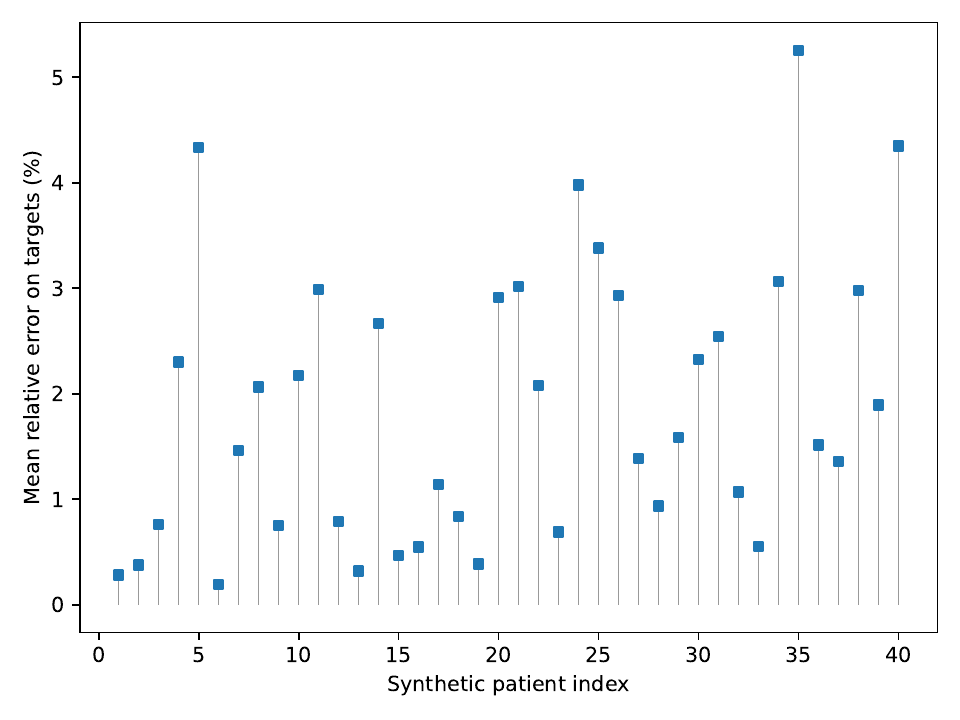}\\
		(c)
	\end{tabular}
	\caption{Resulting mean noise on concentrations for synthetic cohort with an initial noise of order (a) 1\% (b) 3\% and (c) 5\%.} 
	\label{fig:resulting-noise-targets}
\end{figure}

Let $y_p^\star\in\mathbb{R}^5$ denote the reference outlet vector for patient $p$, and $y_p^{\mathrm{pert}}$ the corresponding output obtained with perturbed coefficients. For each component we compute relative errors and summarize them by the mean over the five targets. Averaged over the 40 patient cohort, the mean relative error on the outputs is approximately $0.38\%$ for $1\%$ coefficient perturbations, $1.12\%$ for $3\%$, and $1.87\%$ for $5\%$. These values remain of the same order as the imposed parameter noise and are slightly smaller on average. Even in the worst cases, the maximal mean error across patients remains bounded (below $1\%$, $3.1\%$, and $5.3\%$ respectively). This behavior indicates that, in the considered regime, the forward map does not amplify uncertainties in $(d_{\mathrm{Ca}},d_{\mathrm{Ci}})$ and exhibits a moderate Lipschitz-type sensitivity with respect to these parameters. Figure~\ref{fig:resulting-noise-targets} below exhibits the mean noise on the targets for each synthetic patient.

A finer analysis reveals that the five observable quantities do not respond uniformly to parameter perturbations. Averaged over patients, the largest sensitivities are consistently observed on the citrate-related targets, while albumin-related quantities remain weakly affected. For instance, for $1\%$ perturbations, the mean relative errors are approximately $0.88\%$ and $0.78\%$ for citrate and calcium-citrate on the output, compared to about $0.11\%$, $0.01\%$, and $0.10\%$ for other concentrations. Similar proportions are observed for $3\%$ and $5\%$ perturbations, with citrate and calcium-citrate reaching mean errors of order $2.6\%$-$4.3\%$, while the remaining targets stay below $0.6\%$ on average. This distribution is consistent with the structure of the reaction network: citrate and calcium-citrate species depend more directly on membrane diffusion, whereas albumin-related components are only indirectly affected. Figures~\ref{fig:per-target-sensitivity-1}, \ref{fig:per-target-sensitivity-3} and \ref{fig:per-target-sensitivity-5} below show the noise on the output for each chemical species

\begin{figure}[H]
	\centering
	\begin{tabular}{cc}
		\includegraphics[width=0.5\textwidth]{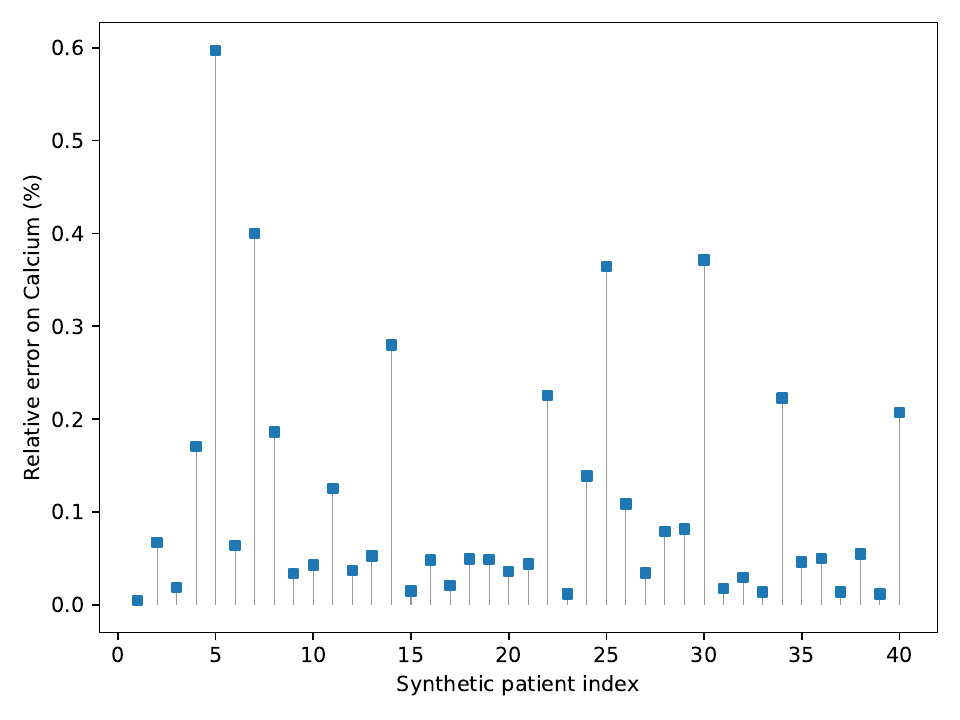}& \includegraphics[width=0.5\textwidth]{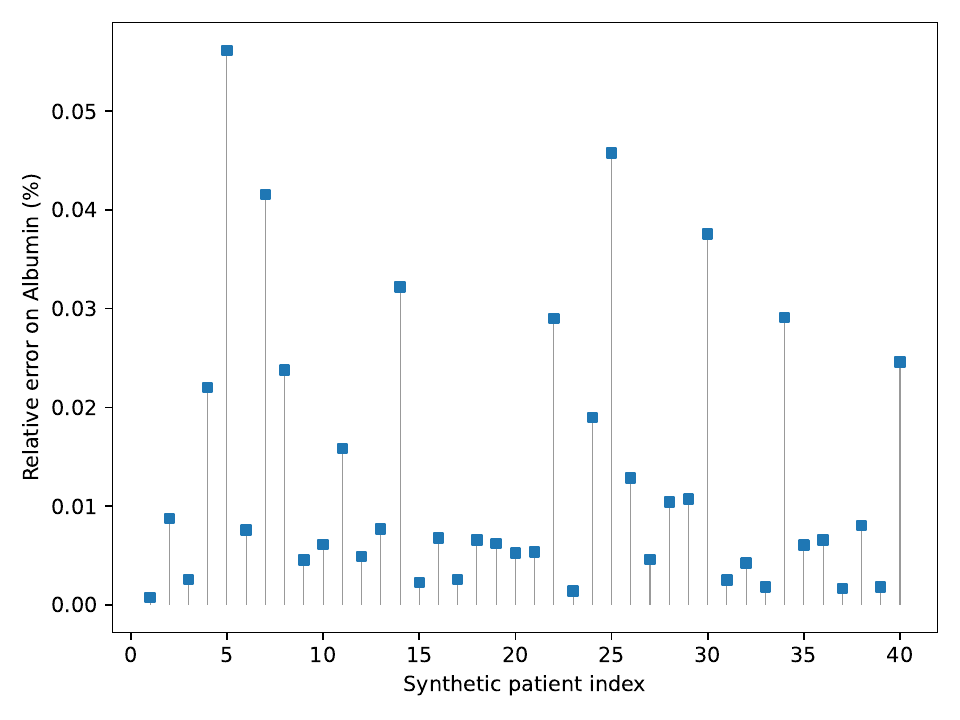}\\
		(a)&(b)
	\end{tabular}\\
	\begin{tabular}{cc}
		\includegraphics[width=0.5\textwidth]{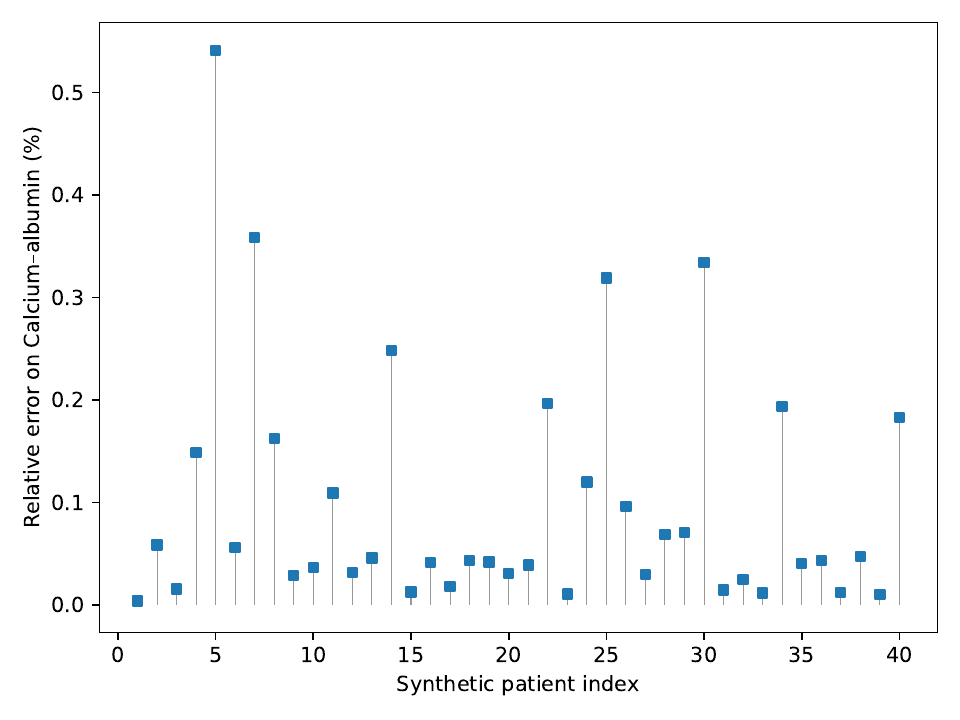}&\includegraphics[width=0.5\textwidth]{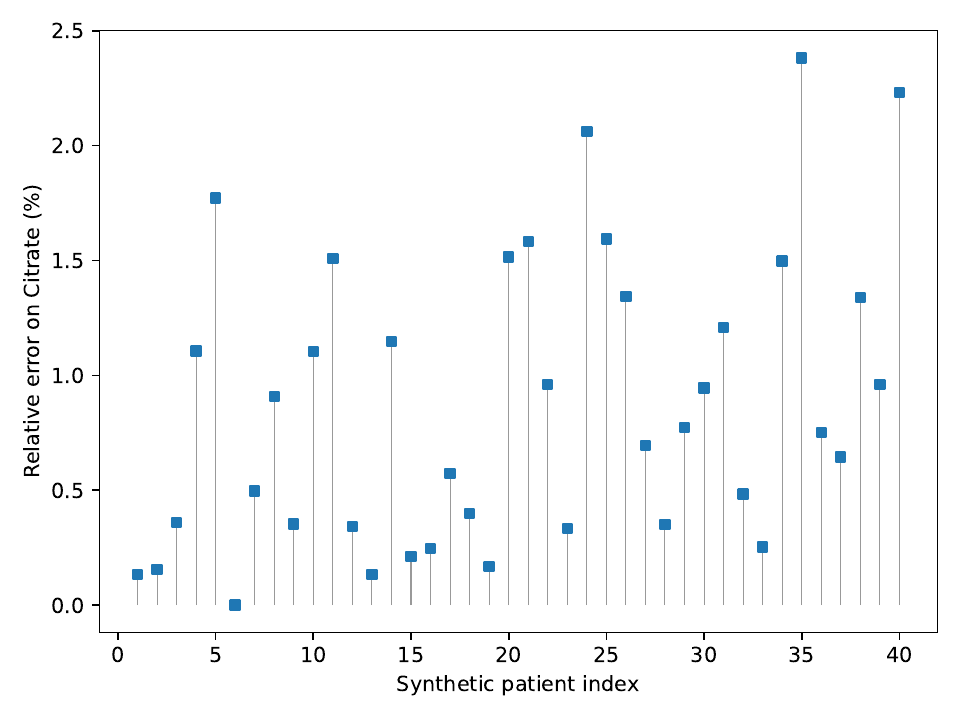}\\
		(c)& (d)
	\end{tabular}\\
	\begin{tabular}{c}
		\includegraphics[width=0.5\textwidth]{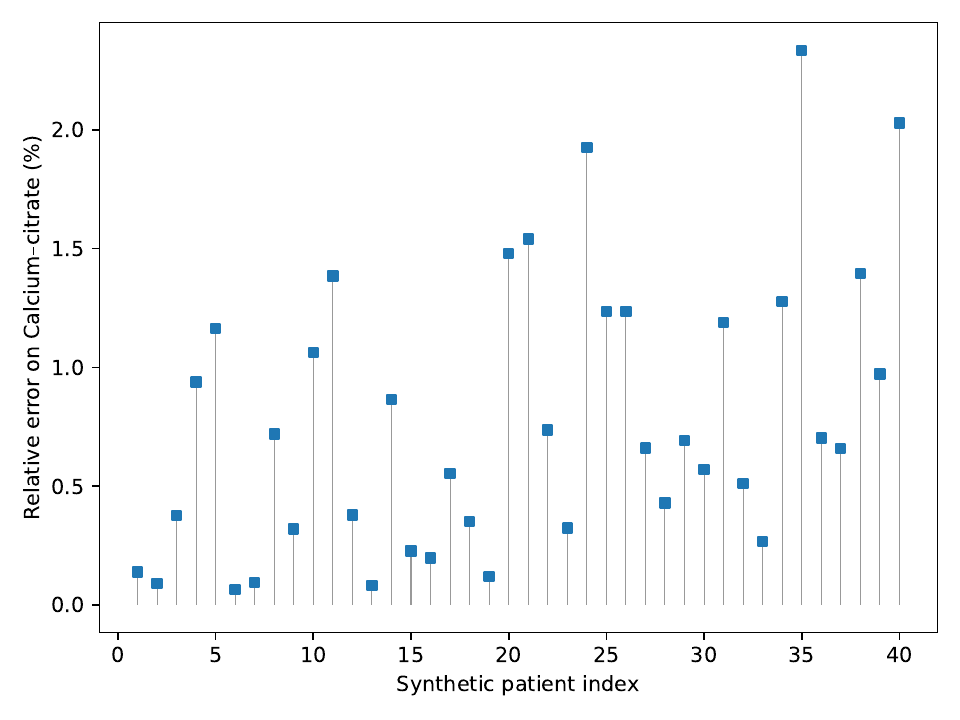}\\
		(e)
	\end{tabular}
	\caption{Resulting noise on each species for 1\%}\label{fig:per-target-sensitivity-1}
\end{figure}

\begin{figure}[H]
	\centering
	\begin{tabular}{cc}
		\includegraphics[width=0.5\textwidth]{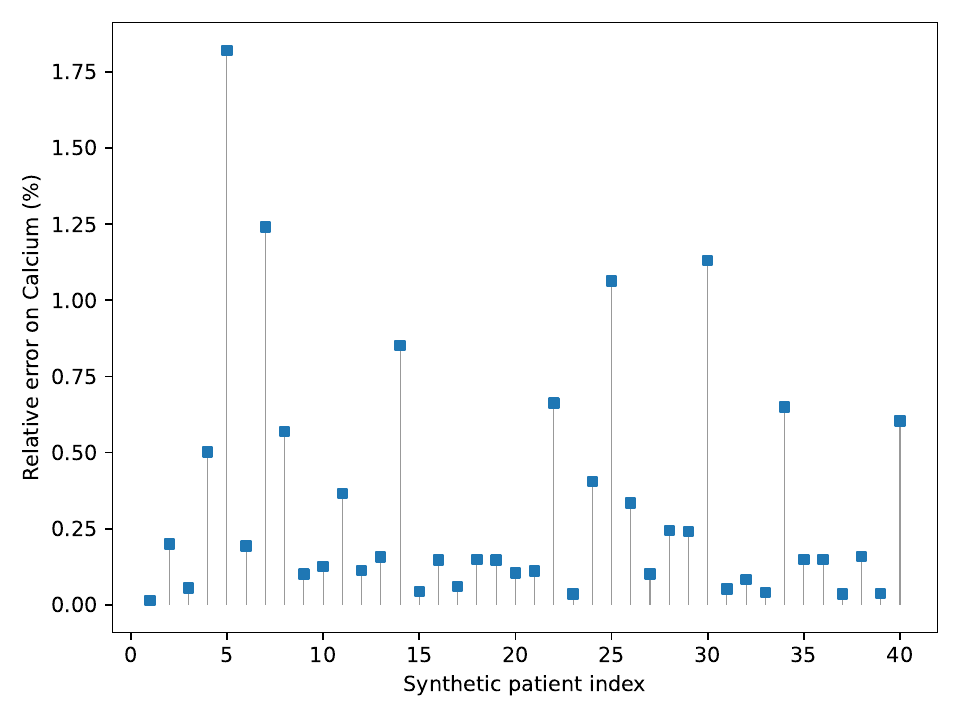}& \includegraphics[width=0.5\textwidth]{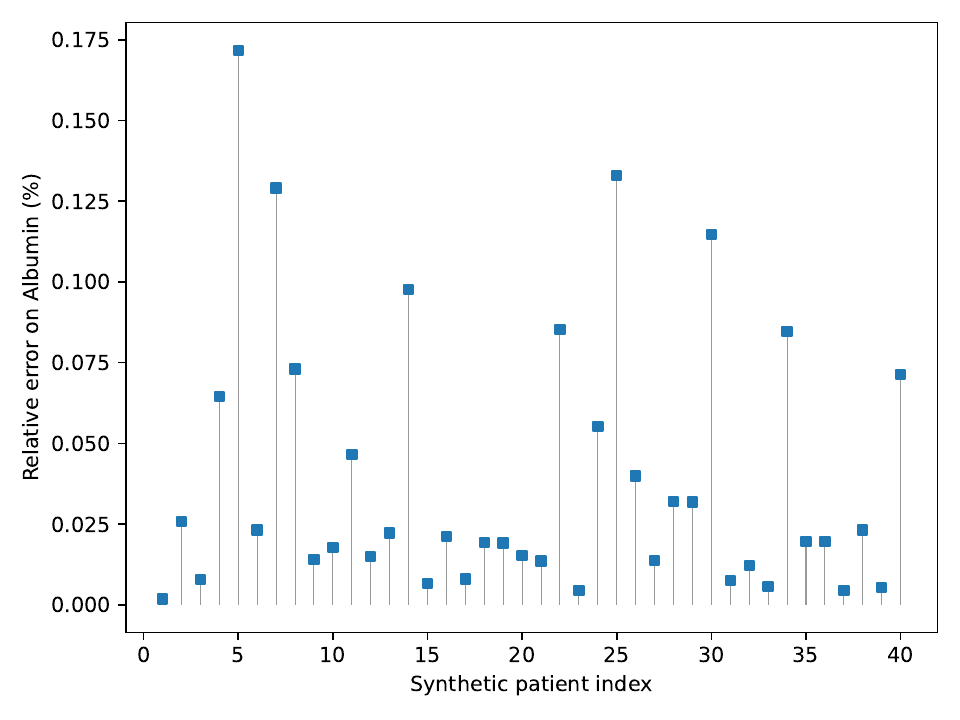}\\
		(a)&(b)
	\end{tabular}\\
	\begin{tabular}{cc}
		\includegraphics[width=0.5\textwidth]{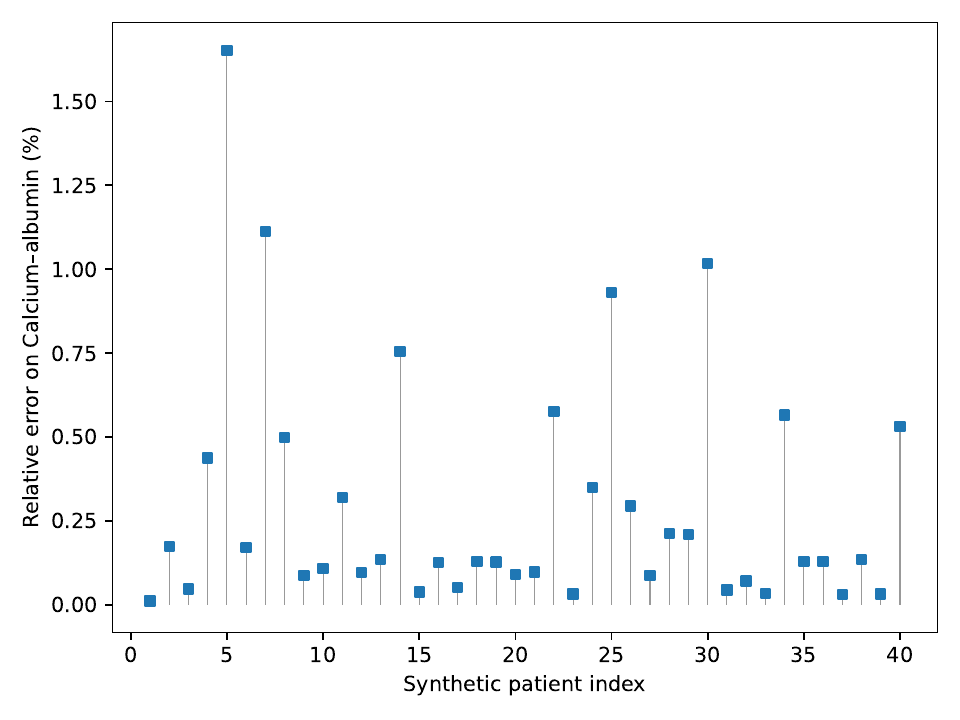}&\includegraphics[width=0.5\textwidth]{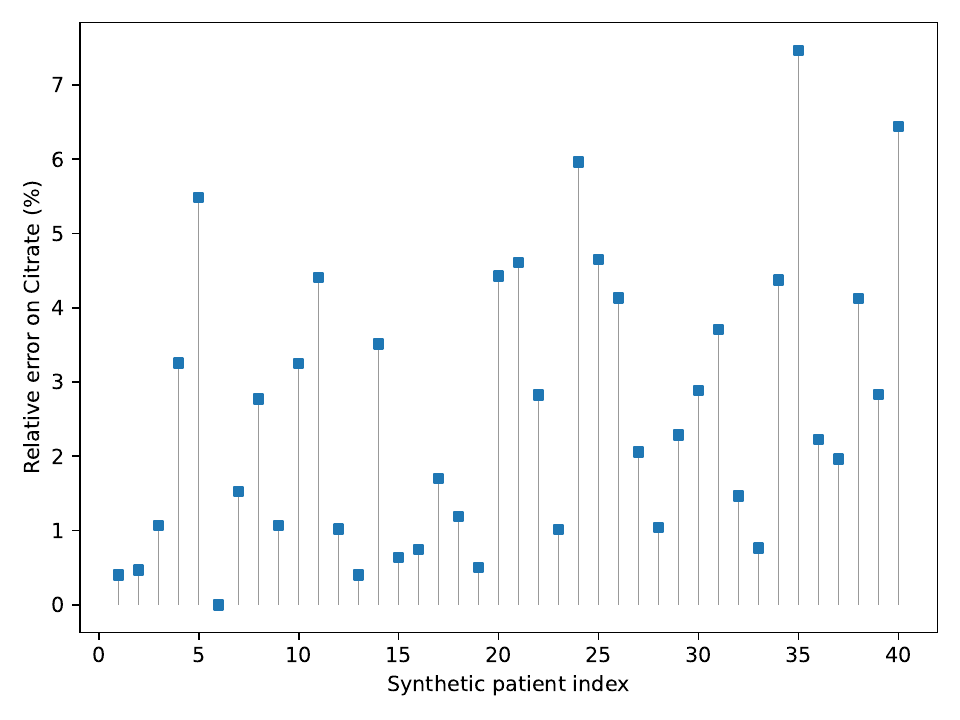}\\
		(c)& (d)
	\end{tabular}\\
	\begin{tabular}{c}
		\includegraphics[width=0.5\textwidth]{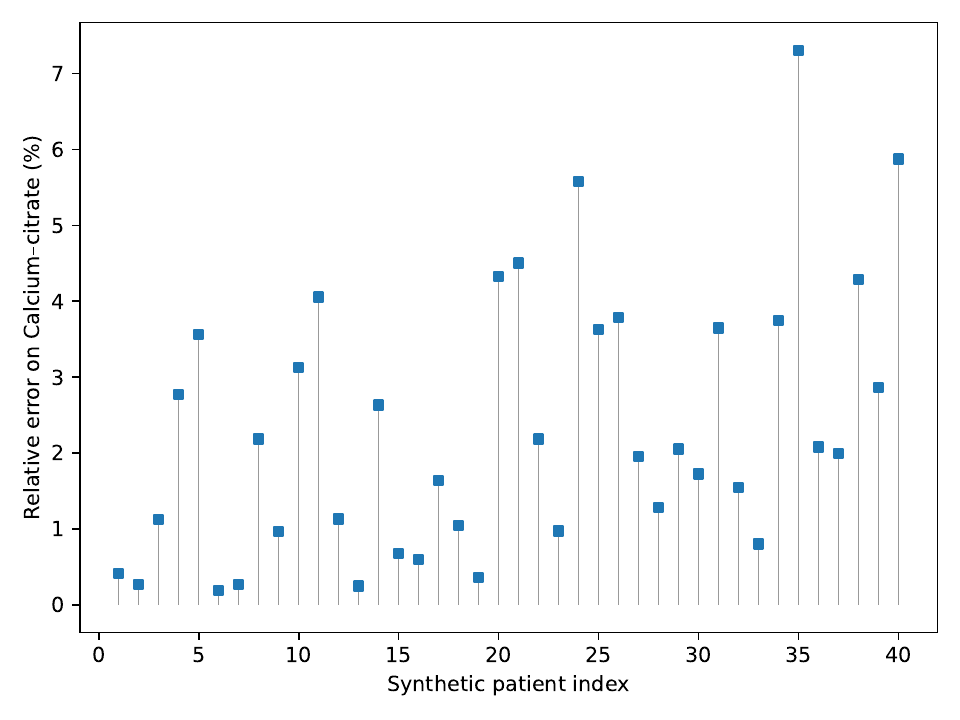}\\
		(e)
	\end{tabular}
	\caption{Resulting noise on each species for 3\%}
	\label{fig:per-target-sensitivity-3}
\end{figure}

\begin{figure}[H]
	\centering
	\begin{tabular}{cc}
		\includegraphics[width=0.5\textwidth]{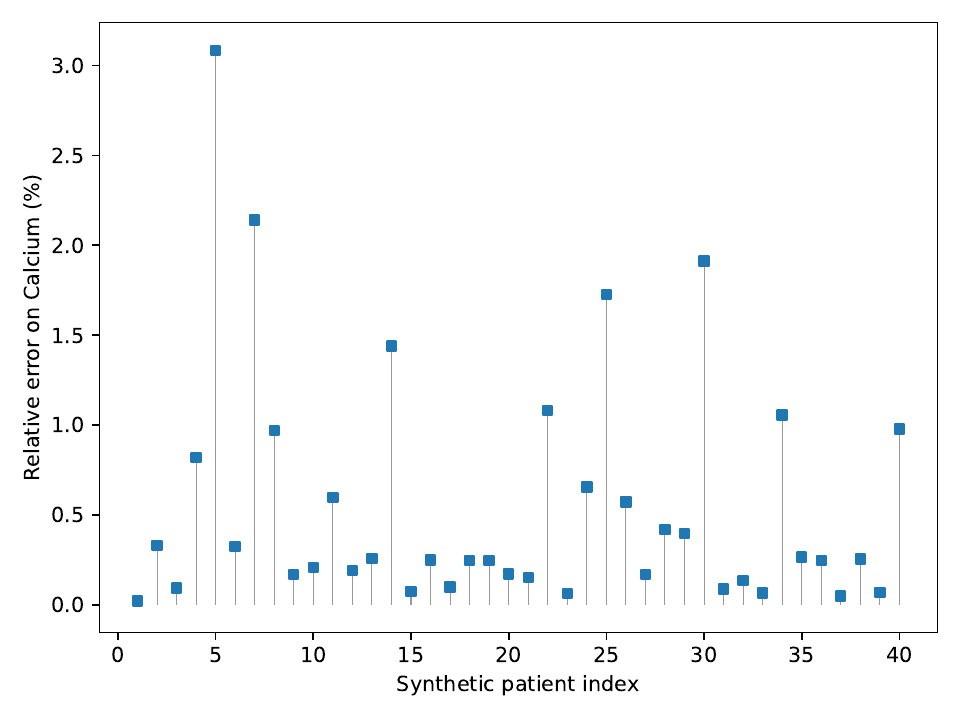}& \includegraphics[width=0.5\textwidth]{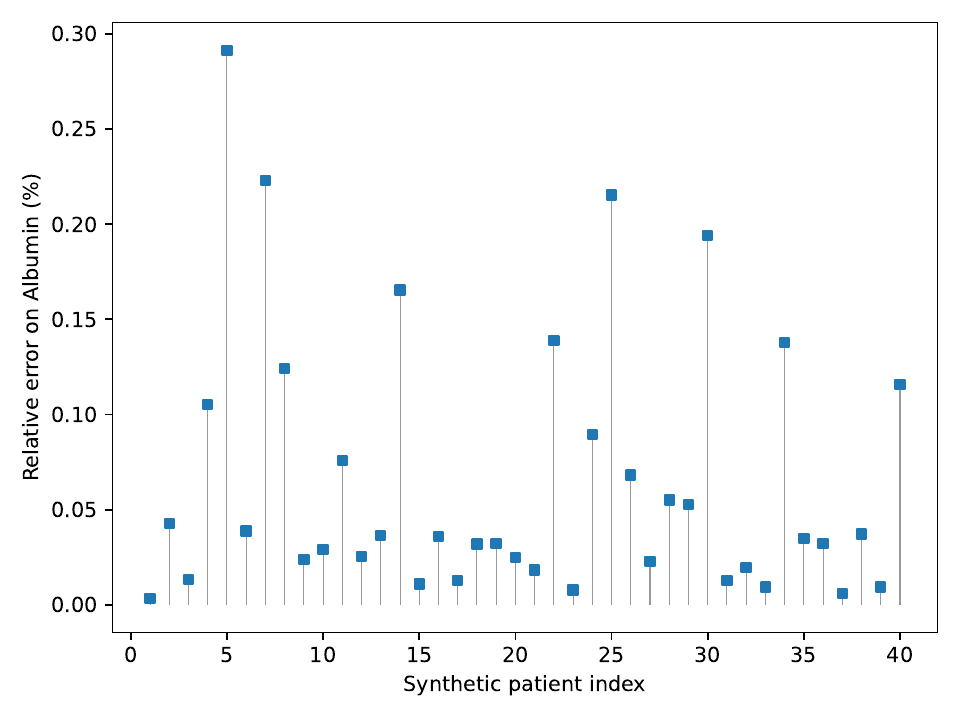}\\
		(a)&(b)
	\end{tabular}\\
	\begin{tabular}{cc}
		\includegraphics[width=0.5\textwidth]{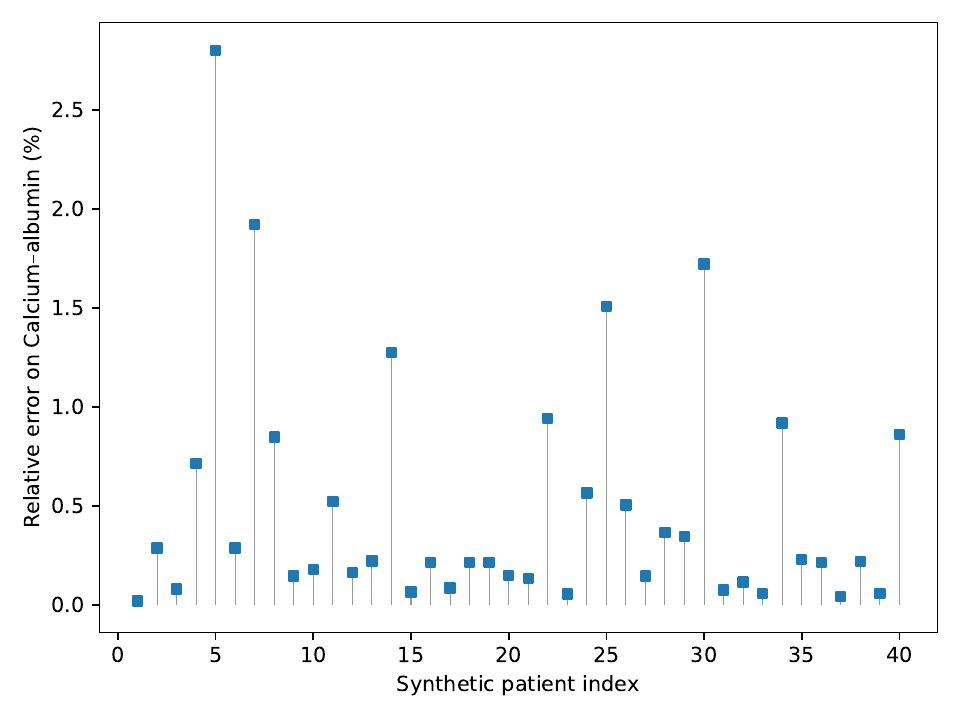}&\includegraphics[width=0.5\textwidth]{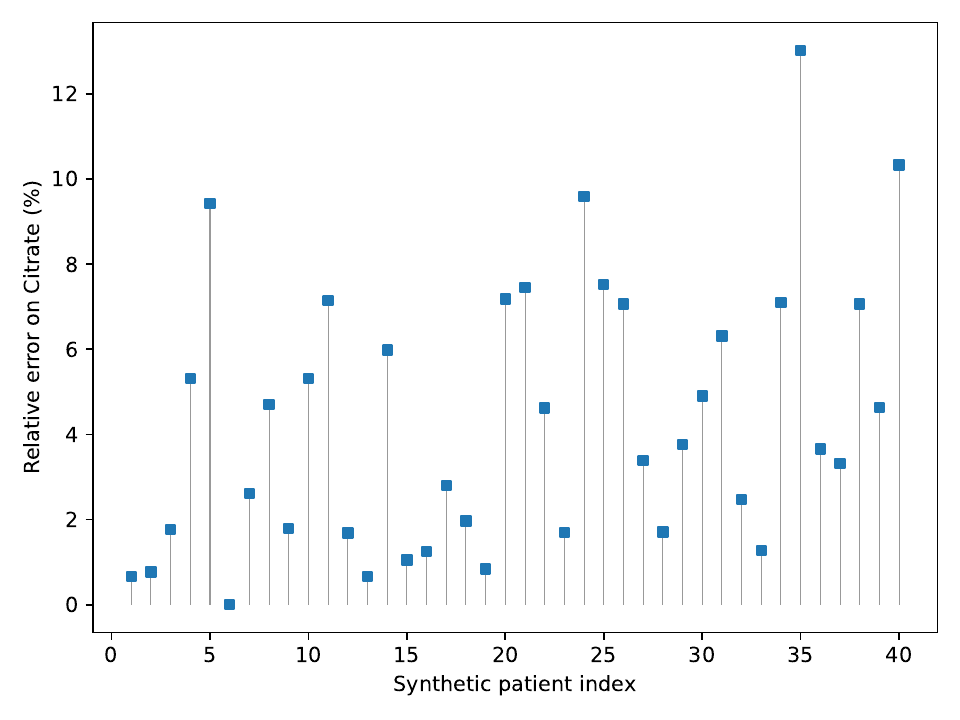}\\
		(c)& (d)
	\end{tabular}\\
	\begin{tabular}{c}
		\includegraphics[width=0.5\textwidth]{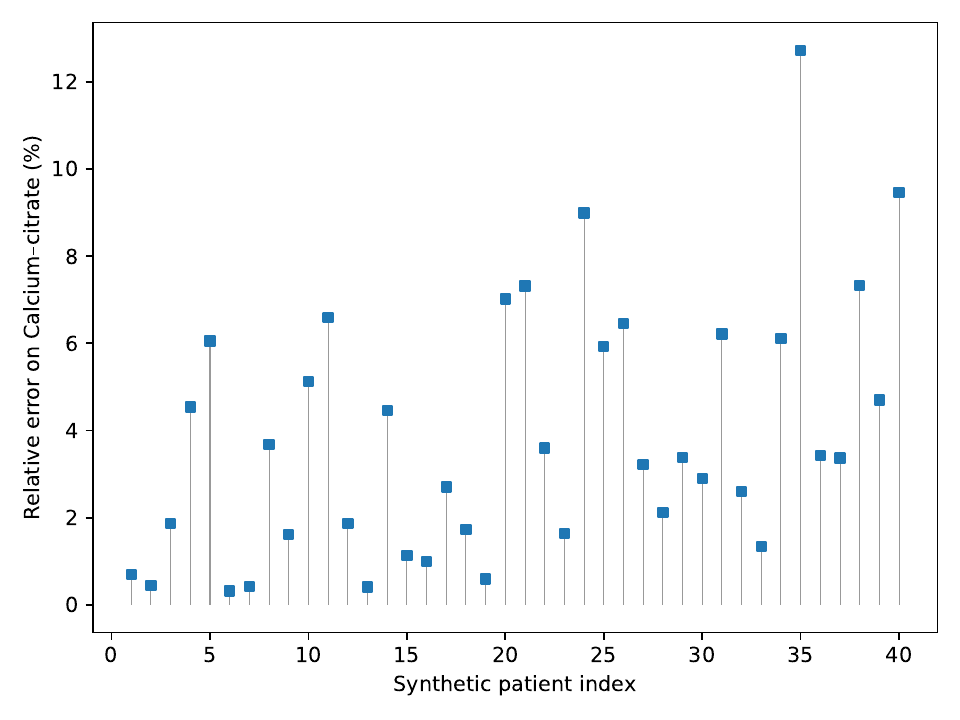}\\
		(e)
	\end{tabular}
	\caption{Resulting noise on each species for 5\%}
	\label{fig:per-target-sensitivity-5}
\end{figure}

Overall, this sensitivity study shows that moderate perturbations of $(d_{\mathrm{Ca}},d_{\mathrm{Ci}})$ induce comparable or smaller relative variations in the model outputs. The absence of amplification supports the numerical stability of the forward problem and reinforces the numerically well-posed condition of the inverse identification in the multi-patient setting. Moreover, the target-wise analysis highlights which observables carry the strongest information about membrane transport properties, providing guidance for weighting strategies in the least-squares functional. Also, the heterogeneous sensitivity observed across the five targets provides a direct explanation for the anisotropic geometry of the objective functional reported in Figure~\ref{fig:identifiability-geometry}. Since citrate-related observables exhibit the strongest dependence on membrane diffusion, variations in $d_{\mathrm{Ci}}$ produce comparatively steep changes in $\mathcal J$, leading to a well-defined minimum in this direction. In contrast, the weaker response of the remaining targets to variations in $d_{\mathrm{Ca}}$ results in a flatter landscape along this axis. This imbalance naturally generates the elongated valley structure observed in the $(d_{\mathrm{Ca}},d_{\mathrm{Ci}})$ plane, with strong identifiability in $d_{\mathrm{Ci}}$ and partial correlation in $d_{\mathrm{Ca}}$. From an inverse-problem perspective, this confirms that the geometry of $\mathcal J$ is primarily driven by the differential sensitivity of the observable species, and that multi-patient aggregation is essential to stabilize the reconstruction along the weakly informed calcium direction.

\section{Application to Clinical Data}\label{sec:real-data-inverse}

We now turn to the application of the proposed inverse identification framework to clinical data obtained in the context of the MARC study (ClinicalTrials.gov identifier: NCT04530175), conducted at CHU Gabriel Montpied (Clermont-Ferrand, France). The cohort consists of twenty-two patients undergoing intermittent hemodialysis, for whom inlet concentrations and flow rates are available, together with outlet blood concentrations of calcium and citrate. As described in the introduction, these data provide sparse but clinically relevant observations of membrane transport processes under real operating conditions. Our objective is to identify effective membrane diffusion coefficients for calcium and citrate, denoted by
\begin{equation*}
	\beta = (d_{\mathrm{Ca}}, d_{\mathrm{Ci}}),
\end{equation*}
by minimizing the multi-patient least-squares functional~\eqref{eq:Jmulti}, which aggregates model-data misfits over all patients and measured chemical species.

\begin{figure}
	\begin{tabular}{cc}
	\includegraphics[width=0.45\textwidth]{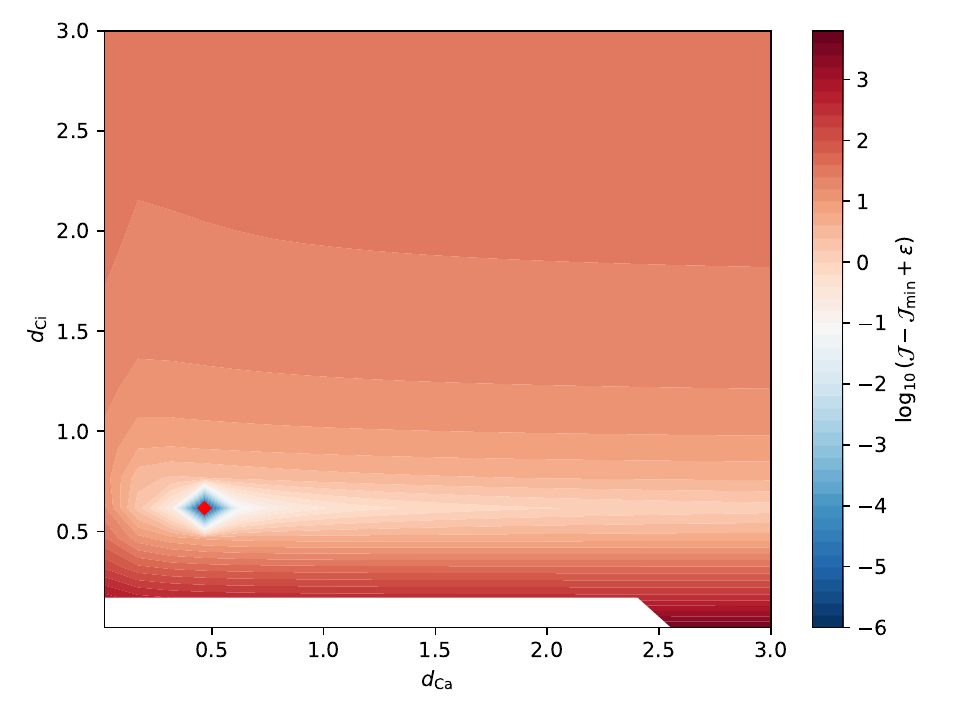} & \includegraphics[width=0.55\textwidth]{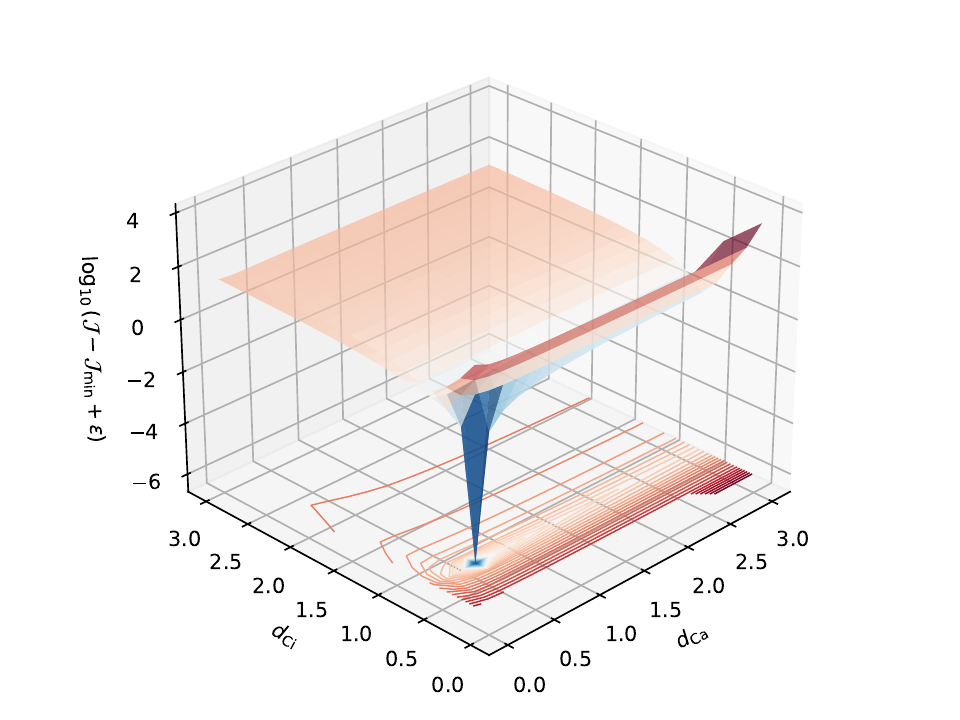}\\
	(a) & (b)
	\end{tabular}
	\caption{Coarse Grid Search in logarithmic scale (a) two-dimensional representation, the grid minimum is the red dot (b) surface representation in 3D.}
	\label{fig:coarse-grid-search}
\end{figure}

As a first exploratory step, we perform a coarse grid search over the  rectangular parameter domain $[0.02,3] \times [0.02,3]$ in the $(d_{\mathrm{Ca}}, d_{\mathrm{Ci}})$ parameters plane. This global scan reveals a pronounced valley-shaped structure of the objective functional, characterized by a relatively sharp sensitivity with respect to $d_{\mathrm{Ci}}$ and a much flatter dependence on $d_{\mathrm{Ca}}$. In particular, the lowest values of $\mathcal{J}$ are consistently observed in a neighborhood of $d_{\mathrm{Ci}} \approx 0.6$ relevant with the previous observations about a strong identifiability of the citrate diffusion coefficient, while calcium diffusion appears less tightly constrained see Figure~\ref{fig:coarse-grid-search}. 

To assess the robustness of this geometric structure and to exclude the presence of local oscillations or secondary minima, we subsequently perform a localized grid search restricted to a narrower region around the valley $[0.25,2.50] \times [0.52, 0.72]$. The resulting refined exploration confirms that the minimum of $\mathcal{J}$ lies within this corridor and that the objective surface is smooth in this region, with no evidence of competing local minima.
% as described in Figure \gl{ref}. 

These observations provide strong numerical support for the stability of the inverse problem in the vicinity of the identified valley. Building upon these grid-based results, we initialize a derivative-free Powell optimization within the valley region. The algorithm converges toward a minimizer
\begin{equation*}
	\beta^\ast = (d_{\mathrm{Ca}}, d_{\mathrm{Ci}}) \approx (0.5043764449030683, 0.6067775575575455),
\end{equation*}
 
which is fully consistent with both the global and localized grid searches. The corresponding value of the objective functional is approximately
\begin{equation*}
	\mathcal{J}(\beta^\ast) \approx 12.279445438717667.
\end{equation*}
While this value represents the lowest misfit attained within the explored parameter space, it remains significantly above zero. This residual discrepancy should not be interpreted as a failure of the identification procedure, but rather reflects several intrinsic limitations of the clinical inverse problem. First, outlet concentration measurements are affected by experimental noise and inter-session variability. Second, the forward model relies on effective transport parameters that lump together membrane microstructure, unresolved physicochemical effects, and modeling simplifications. We note that part of the uncertainty in the estimated parameters may also stem from device-related variability, in the sense that although the study was conducted using hollow-fiber dialyzer cartridges equipped with synthetic polysulfone membranes from the same manufacturer, different cartridge series may have been used in clinical practice, and the effective diffusion properties of the membrane can vary across series. Finally, inter-patient physiological variability cannot be fully captured by a single pair of global coefficients. Consequently, the identified parameters must be understood as effective diffusion coefficients, providing the best compromise in a least-squares sense across the entire cohort. In this context, the relatively elevated value of $\mathcal{J}(\beta^\ast)$ quantifies a combined effect of measurement uncertainty, model mismatch, and biological heterogeneity, rather than purely numerical suboptimality. In this multi-patient setting, deriving confidence intervals on the estimated parameters directly from $\mathcal{J}(\beta^\ast)$ would therefore be misleading, as it would implicitly attribute all discrepancies to observational noise, while significant contributions arise from modeling assumptions and population heterogeneity. As a result, the present framework provides effective population-averaged transport coefficients together with qualitative identifiability information, but does not support a rigorous uncertainty quantification at this stage. Achieving statistically meaningful confidence bounds would require extensive forward evaluations or sampling-based strategies, which are currently computationally prohibitive due to the cost of the underlying partial differential equation model. This limitation motivates the development of surrogate models or reduced-order representations, as discussed in the concluding section.

Despite these limitations, the inverse pipeline exhibits a coherent and reproducible behavior on real clinical data. More precisely, the geometry of the objective functional is well structured, the identified valley is stable across grid resolutions, and the derivative-free optimization consistently converges toward the same parameter region. These features support the practical identifiability of citrate transport and provide a meaningful estimate of calcium diffusion under clinical operating conditions. Beyond parameter calibration, this result constitutes a first step toward predictive patient-specific simulations. Once effective membrane transport properties are identified, the model may be used to explore the impact of alternative dialysate compositions or operating regimes on calcium balance, thereby opening perspectives for personalized dialysis strategies.

\section{Conclusion and perspectives}
\label{sec:conclusions}

In this work, we developed a multi-patient inverse modeling framework for the identification of effective calcium and citrate diffusion coefficients in hollow-fiber hemodialysis devices. Building upon previous modeling efforts, a first methodological contribution is the direct computation of steady-state solutions of the coupled convection-reaction-diffusion system using a Newton method, rather than time-marching transient simulations. This stationary formulation significantly reduces computational cost while preserving accuracy, making repeated forward evaluations feasible within an inverse identification loop. A second major contribution is the extension from single-patient calibration to a multi-patient inverse setting. Using a synthetic cohort constructed from clinical statistics, we formulated a global least-squares objective aggregating residuals across several patients and solved the resulting low-dimensional nonlinear optimization problem using derivative-free methods. Unsurprisingly, numerical experiments with exact synthetic data demonstrate that aggregating information from multiple patients restores numerical identifiability of the diffusion coefficients and mitigates parameter correlation effects that are observed in single-patient inversions. In particular, the combined analysis of grid-based exploration, sensitivity studies, and Powell optimization highlights a pronounced anisotropy in the objective landscape, with strong identifiability in the citrate diffusion direction and weaker sensitivity with respect to calcium.

The robustness study conducted on synthetic cohorts further strengthens the practical relevance of the proposed inverse framework. By introducing controlled multiplicative perturbations on the outlet concentrations, we demonstrated that moderate measurement noise (up to $5\%$ relative standard deviation) does not induce a significant amplification in the reconstructed effective diffusion coefficients. While the dispersion of the estimates increases with the noise level, the identified parameters remain concentrated in a neighborhood of the ground-truth values, with a variability that reflects the intrinsic anisotropy of the objective landscape. In particular, the weaker practical identifiability of $d_{\mathrm{Ca}}$ compared to $d_{\mathrm{Ci}}$ consistently manifests under noisy conditions. These results indicate that the multi-patient formulation provides a stabilizing effect and that the inverse problem remains well-behaved in clinically realistic noise regimes.

Beyond exploring numerical identifiability and robustness to noise on the targets, we investigated the numerical sensitivity of the forward model with respect to perturbations of the diffusion coefficients. By propagating controlled parameter noise through the full forward solver, we showed that relative variations in the observable outputs remain of the same order as the imposed coefficient perturbations, and often slightly smaller. This behavior indicates a moderate conditioning of the forward map and supports the numerical stability of the inverse problem in the considered regime. Target-wise analysis further revealed that citrate-related observables dominate the sensitivity, providing a mechanistic explanation for the valley structure observed in the objective functional.

When applied to real clinical data, the inverse procedure yields effective diffusion coefficients that are consistent across patients and compatible with the geometric structure observed in the synthetic setting. The non-zero residual value of the objective functional reflects several unavoidable factors, as measurement uncertainty, inter-patient physiological variability, and modeling approximations inherent to the reduced partial differential equations description. In particular, the use of effective diffusion parameters can be interpreted as a population-averaged representation of heterogeneous transport phenomena occurring at the fiber scale. Despite these simplifications, the identified coefficients lie in a coherent region of parameter space and reproduce the observed outlet concentrations with satisfactory accuracy. This suggests that the proposed framework captures the dominant transport mechanisms while providing clinically interpretable effective parameters.

From a practical perspective, the main computational bottleneck of the present approach remains the repeated evaluation of the high-fidelity forward model, which involves solving coupled fluid dynamics and nonlinear transport equations for each patient. While tractable in the current two-parameter setting, this cost may become prohibitive for larger cohorts, extended parameter spaces, or real-time clinical applications. A natural perspective is therefore the construction of surrogate models to approximate the parameter-to-observable map. Physics-informed neural networks (PINNs) offer a promising route in this direction, as they can embed the governing equations directly into the training process and provide fast evaluations once trained. In the present context, PINNs could be used either to emulate the full forward solver or to approximate the reduced mapping $(d_{\mathrm{Ca}},d_{\mathrm{Ci}})\mapsto F_p(\beta)$ for each patient, enabling rapid inverse iterations and uncertainty quantification through ensemble evaluations. As an alternative or complementary approach, polynomial or low-rank approximation techniques may also be considered to build reduced-order surrogates of the forward map. Methods based on sparse polynomial expansions or adaptive low-dimensional representations offer strong theoretical guarantees for smooth parametric partial differential equations and can provide interpretable approximations with explicit error control. Such approaches could be particularly attractive in the present low-dimensional parameter setting and would allow systematic exploration of the objective landscape at negligible computational cost once the surrogate is constructed. In particular, Physics-Informed Neural Networks provide a flexible surrogate framework that has already proven effective in several biomedical and microbiological modeling settings, see e.g. the work of P.J.Hossie, B.Laroche, T.Malou, L.Perrin, T.Saigre \& L.Sala~\cite{hossie-laroche-malou-perrin-saigre-sala-2025}.

Overall, the present framework enables the quantitative identification, from available clinical data, of effective transport properties governing chemical exchanges across the dialysis membrane. These coefficients can be interpreted as intrinsic characteristics of a given dialyzer cartridge, capturing in an aggregated manner membrane permeability, microstructure, and unmodeled effects. Once such parameters are estimated with sufficient accuracy, they provide a calibrated mechanistic model of solute transfer through the biomembrane. From a clinical and mathematical perspective, this opens the door to an optimal control formulation of dialysis treatment. More precisely, given patient-specific inlet blood concentrations, the calibrated forward model may be used to determine dialysate compositions that optimize target concentrations at the outlet of the dialyzer. In this setting, membrane transport coefficients play the role of fixed system parameters, while dialysate inputs act as control variables. This perspective naturally leads toward personalized dialysis protocols, in which treatment parameters are adapted to individual patients based on predictive simulations, rather than fixed empirical prescriptions.

{}
\end{document}